\documentclass[a4paper, final]{amsart}

\newcommand{\mytitle}[2][]{%
  \gdef\mylongtitle{\MakeUppercase{#2}}
  \gdef\myshorttitle{\MakeUppercase{#1}}
  \gdef\mypdftitle{#1} %v. \version
}
\mytitle[Hilbert polynomial of length functions]{Hilbert polynomial of length functions}
\newcommand*{\mykeywords}{Length function, Hilbert polynomial, Hilbert series,
algebraic entropy}
\title[\myshorttitle{}]{\mylongtitle%
}

\author[A. Fornasiero]{Antongiulio Fornasiero}
\address{Universit\`a di Firenze}
    
\email{antongiulio.fornasiero@gmail.com} 
\urladdr{https://sites.google.com/site/antongiuliofornasiero/}

\thanks{This work was partially in part by GNSAGA of INdAM and by PRIN 2017,
  project 2017NWTM8R, Mathematical Logic, models, sets, computability.}

\date{13/02/2020}

\usepackage{ifpdf}
\ifpdf
\usepackage[pdftex, linktocpage=true, 
pdfborder={0 0 0}, plainpages=false,%
pdfpagelabels, pdfdisplaydoctitle=true,%
pdfstartview={XYZ null null 1.0},%
colorlinks=false,%
unicode,%
psdextra% additional math macros are supported in bookmarks
]{hyperref}
\hypersetup{%
  pdfauthor={A. Fornasiero},
  pdftitle={\mypdftitle{}},
  pdfkeywords={\mykeywords},
  pdfsubject={Hilbert series and polynomial from length functions}
}
\else
\usepackage[plainpages=false,linktocpage=true,pdfpagelabels, colorlinks=false
]{hyperref}
\fi
\usepackage{booktabs}
\usepackage{multirow}
\usepackage{amsmath, amsthm, amssymb}
\usepackage{mathtools}
\usepackage{xspace}
\usepackage{comment}
\usepackage[neveradjust]{paralist}
\AtBeginDocument{\setdefaultleftmargin{0pt}{}{}{}{}{}}
\setdefaultenum{\upshape(1)}{\upshape i)}{}{}
\usepackage[T1]{fontenc}
\usepackage[backend=biber, style=alphabetic, sorting=nyt,
 backref = true, backrefstyle = three,
 isbn = false, url = false,  
 minalphanames=4,  maxalphanames=5, maxbibnames=99
 ]{biblatex}
\AtEveryBibitem{%
  \clearlist{language}
  \clearfield{urlyear}
  \clearfield{urlmonth}
} 
\newcounter{saveenum}

\newcommand{\vm}{\mathcal V}
\newcommand{\vmc}{\overline{\vm}}
\newcommand{\lambdadim}{\dim_{\lambda}}
\newcommand{\submod}{\leq}
\newcommand{\ideal}{\lhd}
\DeclareMathOperator{\ent}{ent}
\newcommand{\bj}{\bar \jmath}
\newcommand{\btAz}{\tilde{\gr}(A_0; A)}
\newcommand{\btA}{\tilde{\mathcal B}(\overline A)}
\newcommand{\blowtilde}{\tilde{\mathcal B}}
\newcommand{\lcosmall}{$\lambda$-cofinite\xspace}
\newcommand{\ldeg}{$\lambda$-degree\xspace}
\newcommand{\ldim}{$\lambda$-dimension\xspace}
\newcommand{\tildemu}{\tilde \mu}
\newcommand{\muhat}{\hat \mu}
\newcommand{\thetahat}{\hat \theta}
\newcommand{\Aseq}{\mathfrak A}
\newcommand{\Bseq}{\mathfrak B}
\newcommand{\Cseq}{\mathfrak C}
\newcommand{\Mseq}{\mathfrak M}
\DeclareMathOperator{\Rmod}{\mathit{R}-mod}
\DeclareMathOperator{\Smod}{\mathit{S}-mod}
\DeclareMathOperator{\Tmod}{\mathit{T}-mod}
\newcommand{\Rgeqzinf}{\R_{\geq 0} \cup \{\infty\}}
\newcommand{\lS}{\ensuremath{\lambda_S}}
\newcommand{\lR}{\ensuremath{\lambda_R}}
\newcommand{\lT}{\ensuremath{\lambda_T}}
\newcommand{\clLSN}{c\ell^{S,\lambda}_N}
\newcommand{\lSfg}{\lS-f.g\mbox{.}\xspace}
\DeclareMathOperator{\Rees}{\mathcal B}
\newcommand{\ReesN}{\Rees(\overline N)}
\newcommand{\RM}{\Rees(\overline M)}
\newcommand{\RNd}{\Rees_{\bar \delta}(\overline N)}
\DeclareMathOperator{\filt}{Filt}
\DeclareMathOperator{\gr}{Gr}

\DeclareMathOperator{\Ker}{Ker}
\DeclareMathOperator{\Coker}{Coker}
\providecommand{\longto}{\longrightarrow}
\newcommand*{\intro}[1]{\textbf{#1}}
\newcommand*{\Pa}[1]{\bigl( #1 \bigr)}
\newcommand*{\set}[1]{\{#1\}}
\newcommand*{\abs}[1]{\lvert#1\rvert}
\newcommand*{\card}[1]{\lvert#1\rvert}
\newcommand{\N}{\mathbb{N}}
\newcommand{\Z}{\mathbb{Z}}
\newcommand{\R}{\mathbb{R}}
\newcommand{\Q}{\mathbb{Q}}
\newcommand*{\tuple}[1]{\langle #1 \rangle}
\newcommand{\av}{\bar a}

\newcommand{\cv}{\bar c}

\newcommand{\x}{\bar x}
\newcommand{\y}{\bar y}

\def\Ind#1#2{#1\setbox0=\hbox{$#1x$}\kern\wd0\hbox to
  0pt{\hss$#1\mid$\hss}\lower.9\ht0\hbox to 0pt{\hss$#1\smile$\hss}\kern\wd0}

\def\hyph{\nobreakdash-\hspace{0pt}\relax}

\newcommand{\eg}{e.g\mbox{.}\xspace}
\newcommand{\ie}{i.e\mbox{.}\xspace}

\newcommand{\wrt}{w.r.t\mbox{.}\xspace}

\newtheorem{lemma}{Lemma}[section]
\newtheorem{thm}[lemma]{Theorem}
\newtheorem{corollary}[lemma]{Corollary}
\newtheorem{conjecture}[lemma]{Conjecture}
\newtheorem{proposition}[lemma]{Proposition}
\newtheorem{open problem}[lemma]{Open problem}

\newtheorem*{fact*}{Fact}
\newtheorem{fact}[lemma]{Fact}

\theoremstyle{remark}

\newtheorem{claim}{Claim}
\newtheorem*{claim*}{Claim}
\newtheorem{exercise}[lemma]{Exercise}

\theoremstyle{definition}
\newtheorem{definition}[lemma]{Definition}

\newtheorem{remark}[lemma]{Remark}
\newtheorem{final remark}[lemma]{Final remark}
\newtheorem{example}[lemma]{Example}
\newtheorem{examples}[lemma]{Examples}

\newenvironment{adef}[1]{%
  \par\noindent\textbullet\ 
}{%
}

\newenvironment{sentence}[1][]{%
  \begin{list}{}{%
    \setlength\topsep{0.5ex}%
    \setlength\leftmargin{\parindent}%
  }%
  \item[#1]
 }
 {\end{list}}

\begin{document}

\begin{abstract}
Let $\lambda$ be a general length function for modules over a Noetherian ring R.
We use $\lambda$ to introduce  Hilbert series and polynomials  for R[X]-modules,
%We study the notion of algebraic entropy for $R[X]$-modules.
%with respect to a length function $\lambda$ on R-modules. 
%We introduce a Hilbert series and a Hilbert polynomial 
measuring the growth rate of~$\lambda$. 
We show that the leading term $\mu$ of the Hilbert polynomial is an
invariant of the module, which refines both the algebraic entropy and the
receptive algebraic entropy; its degree is a suitable notion of dimension for
$R[X]$-modules.

% We  also consider various modifications of the Hilbert polynomial, such as the
% intrinsic and the multivariate version of the Hilbert polynomial, and relate
% them to different versions of entropy. 

% Let $\lambda$ be a general length function for modules over a Noetherian ring R.
% We use $\lambda$ to define  Hilbert series and polynomials  for R[X]-modules.
% The leading term $\mu$ of any such polynomial is an
% invariant of R[X]-modules, which refines both algebraic entropy and the more
% recently introduced receptive algebraic entropy, and the degree of $\mu$ is a suitable notion of dimension for $R[X]$-modules.

Similar to algebraic entropy, $\mu$ in general is not additive for exact sequence of 
$R[X]$-modules: we demonstrate how to adapt of certain entropy constructions  to this new
invariant.

We also consider multi-variate versions of the Hilbert polynomial.

% Let L be a general length function for modules over a Noetherian ring R.
% We use L to define  Hilbert series and polynomials 
% for R[X]-modules.
% The leading term of any such polynomial is an
% invariant of R[X]-modules, which refines the algebraic entropy. 
\end{abstract}

\maketitle

{\small
\setcounter{tocdepth}{1}
%% only sections in the TOC
\tableofcontents
}

\makeatletter
\renewcommand\@makefnmark%
   {\normalfont(\@textsuperscript{\normalfont\@thefnmark})}
\renewcommand\@makefntext[1]%
   {\noindent\makebox[1.8em][r]{\@makefnmark\ }#1}
% change style of footnotes: (1) instead of 1 
% done here instead than in the preamble to avoid messing with amsart formatting of keywords et similia
\makeatother

\section{Introduction}

Let $R$ be a  commutative ring with unity; for this introduction and
most of the article $R$ will be Noetherian.

A generalized \intro{length} function on the category $\Rmod$ of $R$-modules is a function
\[
\lambda: \Rmod \to \Rgeqzinf
\]
satisfying the following conditions:
\begin{enumerate}
\item $\lambda(0) = 0$;
\item $\lambda(M) = \lambda(M')$ when $M$ and $M'$ are isomorphic;
\item for every exact sequence
$0 \to A \to B \to C \to 0$,
\[
\lambda(B) = \lambda(A) + \lambda(C)
\]
(with the usual rule that $x + \infty = x$);
\item for every $M \in \Rmod$,
\[
\lambda(M) = \sup \set{\lambda(M'): M' \submod M\text{ finitely generated $R$-submodule}}.
\]
\end{enumerate}
Generalized length functions were introduced in \cite{NR:65} and further
studied (among other sources) in \cites{Vamos:68, SVV}: see also
\S\ref{sec:length} for a brief introduction; we will simply say ``length'' or ``lenght function'' in the rest of the paper.

%\footnote{Some of the results may generalize to non-commutative rings}

Fix $1 \leq  k \in \N$ and let $S := R[x_{1}, \dotsc, x_{k}]$. 
Let $M$ be an $S$-module.

% Algebraic entropy is a measure of the complexity of a $S$-modules, which is
% inspired to the classical notions of entropy for dynamical systems. 
% In its original version, it captures the asymptotic growth rate of the length of
% the images of an $R$-module under a linear transformation; later it has been
% generalized in various ways, in particular for $S$-modules (equivalently, for
% modules with $k$ commuting linear transformations) using a generalized length
% function (see \S\ref{sec:entropy}).

% In this paper, we investigate the properties of
% algebraic entropy for $S$-modules;
% we introduce a new tool,  the Hilbert series and the Hilbert polynomial of~$λ$,
% which encodes the information about the length of the successive images of a module. 
% The leading term of the Hilbert polynomial is an invariant of the module, which
% refines both the algebraic entropy and the receptive algebraic entropy. 

We generalize the theory of Hilbert series and Hilbert polynomial for
$S$-modules when $R$ is a field and
the linear dimension is the length function 
(see e.g.\ \cites{Eisenbud,MS:05}),
to the case of $S$-modules with an arbitrary length function~$\lambda$. 
We begin by assuming that $M$ is an $\lambda_{S}$-small module, \ie,  
there exists a finitely generated $R$-submodule $V$ of $M$ with finite $\lambda$-length
such that $SV = M$ (we say that $V$ witnesses that $M$ is $\lS$-small). 
We denote by $S_{n}$ the set of polynomials in $S$ of total
degree less or equal to~$n$,
and consider the formal power series
\[
 \sum_{n = 0}^{\infty} \lambda(S_{n} V) t^{n}
\]
and prove that it is a rational function (of~$t$).
We also show that for large enough $n$, the function $n \mapsto \lambda(S_{n} V)$
is a polynomial, whose leading term $\mu(M)$ is independent of the choice of the
witness~$V$
(Theorem~\ref{thm:Hilbert-poly}): thus, 
$\mu_{\lambda}(M)$ is an invariant of~$M$ that measures the asymptotic growth of $\lambda$
on~$M$, and
refines both the algebraic entropy and the receptive algebraic entropy. 
Moreover, the degree of $\mu(M)$ gives a well-behaved notion of dimension (\wrt $\lambda$) for $S$-modules.

% Under some  assumptions on $M$, we 
% define a \intro{Hilbert series} and a \intro{Hilbert polynomial} for $M$, 
% generalizing the known theory when $R$ is a field and $\lambda$ is equal to the
% linear dimension (see e.g.\ \cites{Eisenbud,MS:05}).
% More precisely, we fix a finitely generated $R$-submodule $V \submod M$ such that $\lambda(V) < \infty$ and $SV
% = M$ (we say that $V$ witnesses that $M$ is $\lS$-small);  
% we  define $S_{n}$ as the set of polynomials in $S$ of total
% degree less or equal to $n$; we prove that  the formal power series
% \[
%  \sum_{n = 0}^{\infty} \lambda(S_{n} V) t^{n}
% \]
%  is a rational function (of $t$), and the
% function
% $n \mapsto \lambda(S_{n} V)$ is equal to a polynomial (for $n \in \N$ large enough).
% The leading
% term $\mu_{\lambda}(M)$ of this polynomial does not depend on the choice of the witness~$V$
% (Theorem~\ref{thm:Hilbert-poly}): thus, 
% $\mu_{\lambda}(M)$ is an invariant of~$M$, and describes the growth rate of $\lambda$ on~$M$.  

\smallskip

In~\S\ref{sec:filter} we review some basic notions and results about graded and
filtered modules over a Noetherian ring. 
In~\S\ref{sec:Hilbert}, we construct the Hilbert series for two classes of
modules: 
graded modules (Theorem~\ref{thm:Hilbert-graded})
 and upward filtered modules (Corollary~\ref{cor:Hilbert-filtered}). 
We choose to work with upward filtered modules instead of the more common
downward filtered modules, because they are more suitable for the applications
in 
\S\ref{sec:growth} and  \S\ref{sec:additivity}
(see also\cite[\S1.3]{KLMP}).

% After  some preliminary definitions and facts about
% graded and filtered modules in~\S\ref{sec:filter},
% in \S\ref{sec:Hilbert} we show the existence of a Hilbert series for  
% the case when $M$ is an $\N$\hyph{}graded module
% (Theorem~\ref{thm:Hilbert-graded}); then we move to the case when $M$ is a
% upward filtered module (Corollary~\ref{cor:Hilbert-filtered}).
% Often Hilbert series are defined for downward filtered modules, but,
% in light of the applications in \S\ref{sec:growth} and  \S\ref{sec:additivity}, it is
% more appropriate for the present purposes to consider upward filtered modules (\cf \cite{KLMP}*{\S1.3}).

% In \S\ref{sec:Hilbert-p} we deduce the existence of a Hilbert polynomial 
% for $\lS$-small modules,
% and show that its leading term is an invatiant of~$M$.

% In \S\ref{sec:general} we generalize the definition of $\mu(M)$ to the case when 
% $M$ is \emph{not} $\lS$-small,
% and show that, under some suitable assumptions on the modules (\ie, that they
% are locally $\lR$-finite: see Def.~\ref{def:lambda-finite}),
% the function $\mu$ is
% additive (Theorem~\ref{thm:additivity}).

In \S\ref{sec:Hilbert-p}, we prove the existence of a Hilbert polynomial for
$\lS$-small modules, and show that its leading term, denoted by $\mu(M)$, is an
invariant of~$M$  

In \S\ref{sec:general},
we extend the definition of $\mu(M)$ to the case when $M$ is \emph{not} $\lS$-small, and
show that  $\mu$ is an additive function on the class of modules that are locally
$\lR$-finite (see Def.~\ref{def:lambda-finite} and Theorem~\ref{thm:additivity}).

The coefficient of the $k$-term of the Hilbert polynomial is (up to a constant factor)
the algebraic entropy of the action of $\N^{k}$ on $M$ 
(see \S\ref{sec:length} for the definition of algebraic
entropy and its main properties). 
Therefore, the additivity of $\mu$ is a refinement of the
known additivity of algebraic entropy 
(see Fact~\ref{fact:add-entropy} and \cites{SVV,SV:15,DFG}).
However, the additivity of algebraic entropy has already been proved under weaker
assumptions: one of the most general results considers the case when the acting
monoid $\N^{k}$ is replaced by a cancellative and amenable monoid 
(and $M$ is locally $\lambda_{R}$-finite): %(as  in Definition~\ref{def:lambda-finite})
see \cites{Virili:19, DFG}. 

% The coefficient of the $k$-term of the Hilbert polynomial is (up to a constant
% factor) the algebraic
% entropy of the action of $\N^{k}$ on $M$ (see \S\ref{sec:length} for the
% definition of algebraic entropy and its main properties).
% Therefore, the additivity of~$\mu$ is a refinement of the known   additivity
% of algebraic  entropy (see Fact~\ref{fact:add-entropy} and \cites{SVV,SV:15,DFG}).
% However, the additivity of algebraic entropy
% has already been proved under weaker assumptions: one of the most general results
% considers the case when the acting monoid $\N^{k}$ is replaced by a cancellative
% and amenable monoid, and 
% $M$ is locally $\lambda_{R}$-finite (as  in Definition~\ref{def:lambda-finite}): see
% \cites{Virili:19, DFG}. 

In \S\ref{sec:homogeneous} we show how the usual construction of Hilbert-Samuel
polynomial can be extended to length functions, 
thus obtaining another invariant of~$M$.

In \S\ref{sec:algebra}   we replace $S$ with a
finitely generated $R$-algebra $T$
and define a corresponding Hilbert polynomial for each
$T$-module~$M$: its degree will be an
invariant of~$M$ (while the leading coefficient will depend on the choice of a
set of generators for~$T$).

In \S\ref{sec:receptive} we introduce the $d$-dimensional entropy
%(essentially, the $d$th coefficient of the Hilbert polynomial)
as a generalization of the receptive entropy in \cite{BDGS1}, and relate it
to the Hilbert polynomial.

\smallskip

It is well known that the algebraic entropy of modules that are not
$\lR$-finite may fail to be additive, which is a desirable property. 
To overcome this limitation, some alternative notions of algebraic entropy have
been introduced in the literature. We will explore how similar adaptations can
be applied to~$\mu$, the leading term of the Hilbert polynomial. 

In \S\ref{sec:hat}, we use a technique from \cite{Vamos:68} to define $\muhat$,
an additive function on all $S$-modules that extends~$\mu$ on locally $\lR$-finite
modules. 
This construction also works for the ($d$-dimensional) entropy. 

In \S\ref{sec:intrinsic}, we define the ``intrinsic'' Hilbert polynomial, which
is related to the intrinsic algebraic entropy introduced in~\cite{DGSV:15}: see
\S\ref{sec:entropy-intrinsic}. 
We obtain another invariant $\tildemu$ from the intrinsic Hilbert polynomial. 
We conjecture that $\tildemu$ is additive on $S$-modules, and prove that it is
sub-additive. 
Under this conjecture, we have two additive invariants, $\muhat$ and $\tildemu$,
which may differ in general. 

Moreover, in \S\ref{sec:fine}, we consider a finer version of the Hilbert series where the grading is given by a suitable monoid $\Gamma$ instead of $\N$.

% It is well known that on modules that are not $\lambda_{R}$-finite the algebraic
% entropy might not be additive.
% To ameliorate this ``defect'' some variants of the algebraic entropy have been
% proposed: we will show that similar ideas can be employed also for~$\mu$.

% In \S\ref{sec:hat} we employ a technique in \cite{Vamos:68} and
% use $\mu$ to define an invariant $\muhat$ which is a
% length function on all $S$-modules: %(assuming that $R$ is Noetherian): 
% a similar
% construction works for the ($d$-dimensional) entropy.

% In \S\ref{sec:intrinsic} we define the ``intrinsic'' Hilbert polynomial: it
% bears the same relationship to the Hilbert polynomial that the intrinsic
% algebraic entropy (introduced in \cites{DGSV:15}: see \S\ref{sec:entropy-intrinsic}) bears to the ``usual'' algebraic
% entropy.
% From it we extract an invariant $\tildemu$.
% We conjectuere that $\tildemu$ is an
% additive function on $S$-modules, and prove that is sub-additive. 
% Thus, under the above conjecture, we have
% we have  two  additive functions on $S$-modules, $\muhat$
% and $\tildemu$, which in general are different.

% Besides, we consider in \S\ref{sec:fine} a finer version of the Hilbert series 
% where the grading is no longer in $\N$ but in some suitable monoid $\Gamma$.
%Theorem~\ref{thm:Hilbert-orbit-monoid} seems to be new even in the classical case
%when $R$ is a field (and $\lambda$ is the linear dimension).

% replacing conditions like ``$N$ is a finitely generated $R$-module'' with
% ``$N$ has a submodule $N'$ such that $N'$ is finitely generated, and $\lambda(N/N')
% = 0$''

\medskip

We mostly adapt well-known results about Hilbert series and polynomials to our
setting, or prove them by simple arguments. Therefore, we omit some proofs for
brevity.

\smallskip

In most of the results, we  assume that $R$ is Noetherian ring (or at least that the
relevant modules are Noetherian).
%We assume throughout the paper, except in the appendices, that R is a
%Noetherian ring. 

We conjecture that some of our results can be extended to non-Noetherian rings.
%(see the Appendices for some partial results). 
However, we lack a satisfactory notion of ``Noetherianity with respect to
$\lambda$'' for rings, which prevents us from pursuing this direction further. We
remark that, without any Noetherian assumption, a Hilbert polynomial may not
exist (see \S\ref{subsec:ex-noetherian}). 
On the other hand, the algebraic entropy and its
intrinsic version have been studied for non-Noetherian rings (see \eg \cites{SV:15, SV:19}).

\smallskip

We also leave as an open problem the case when either $R$ or $S$ are
non\hyph{}commutative rings. 
See \cite{Northcott:68} for the case
when $R$ is not commutative, and \cites{DFG, Virili:19} for the case when $S$
is not commutative.
Some partial results are in Appendix~\ref{sec:noncommutative-hat}. 

We believe that most results (except possibly Proposition~\ref{prop:cosmall} and its
corollaries) can be generalized to the case when $R$ is non-commutative, but we do
not explore this possibility here due to our limited expertise in
non-commutative rings.

% \medskip
% Most of the proof are either adaption of ``classical'' proofs about Hilbert
% series and polynomials to this setting, or quite straightforward: hence we
% sometimes skip them.

% \smallskip

% In most of the results, we  assume that $R$ is  Noetherian (or at least that the
% relevant modules are Noetherian).
% We think it may be possible to generalize the results to some non Noetherian rings
% (see the Appendices); however, we don't have (yet) a working
% theory of ``Noetherianity-up-to-$\lambda$'' for rings, so we could not go far in this
% direction of inquiry.
% Without any assumption of Noetherianity  a Hilbert polynomial may not exist (see
% \S\ref{subsec:ex-noetherian}): on the other hand, the algebraic entropy (and its
% intrinsic version) for non-Noetherian
% rings has already been studied (see \eg \cites{SV:15, SV:19}).

% It would also be interesting to consider the cases when either the ring $R$
% or the algebra $S$ are not commutative  (\cf \cite{Northcott:68} for the case
% when $R$ is not commutative, and \cites{DFG, Virili:19} for the case when $S$
% is not commutative): some partial results are in
% Appendix~\ref{sec:noncommutative-hat}; we think that most results (except
% possibly Proposition~\ref{prop:cosmall} and its corollaries) generalize to the case when $R$ is
% non-commutative; however, we mostly avoided this line of inquiry due to our lack
% of familiarity with non-commutative rings.

\subsection*{Acknowledgments}
We thank Simone Virili and Giorgio Ottaviani for the useful discussions,
and the anonymous referee for the care and great amount of work that he
put in refereeing this article.

\section{Preliminaries, assumptions, and notation}

\subsection{Notation}
$\N$ denotes the set of natural numbers, including $0$.\\
$\infty$  denotes some element that is greater than any real number.\\
$R$ is a ring (commutative with $1$) and $\lambda$ is a length function on $\Rmod$.\\
We write $I \ideal R$ if $I$ is an ideal of~$R$,
and $A \submod M$ if $A$ is a submodule of the module~$M$.\\
We fix $1 \leq k \in \N$ and denote $\x \coloneqq\tuple{x_{1}, \dotsc, x_{k}}$
 and $S \coloneqq R[\x]$.\\
Given $n \in \N$, $S^{(n)}$ denotes the set of homogeneous polynomials in $S$ of
degree exactly $n$ (plus the $0$ polynomial), while $S_{n}$ is the set of
polynomials in $S$ of degree at most~$n$ (they are both finitely generated
$R$-modules: notice that $S_{0} = R$).%
%\footnote{$S^{(n)}$ and $S^{\leq n}$ might be better notation?}
\\
From \S\ref{sec:Hilbert-p} to the appendices $R$ will be a Noetherian ring
(therefore, $S$ and $S[y]$ will also be Noetherian rings).

\subsection{Polynomial coefficients for some rational functions}
In this subsection we gather some results:
%We will need the following results: 
 probably they are well known, but we could not find a reference.

We fix some $\ell \in \N$ and denote $\bar t \coloneqq \tuple{t_{1}, \dotsc, t_{\ell}}$.

The case $\ell = 1$ of the next proposition is well-known, and the one we will use
for most of the article (see \eg \cite[Ch.11]{AM}).
\begin{proposition}\label{prop:poly-growth}
Let $K$ be a ring  of characteristic~$0$, and
$p(\bar t) \in K[\bar t]$.
Let $\gamma_{1}, \dotsc, \gamma_{\ell} \in \N$.
Define
\[
f(\bar t) := \frac{p(\bar t)}{\prod_{i=1}^{\ell}(1 - t_{i})^{\gamma_{i}}} \in K(\bar t)
\]
and expand $f$ as
\[
f(\bar t) = \sum_{\bar n \in \N^{\ell}} a_{\bar n} \bar t^{\, \bar n} \in K[[\bar t]].
\]
Then, there exists a polynomial $q(\bar t) \in K[\bar t]$ such that:
\begin{enumerate}
\item For every $\bar n \in \N^{\ell}$ large enough
\[
a_{\bar n} = q(\bar n);
\]
\item for every $i= 1, \dotsc, \ell$
\[
\deg_{t_{i}}(q) \le \gamma_{i} - 1
\]
with $\deg_{t_{i}}(q) = \gamma_{i} - 1$
if $K$ is an integral domain and $p \neq 0$.
%(where we set $\deg(0) \coloneqq - 1$).
\end{enumerate}
\end{proposition}

\begin{proof}
It is clear that it suffices to treat the case when $p=1$.
We proceed by induction on $\ell$.
%$$k := \sum_{i=1}^{\ell} \gamma_{i}$.
If $\ell = 0$, then $f=1$, and $q= 0$.
If $\ell = 1$, the result is easy: by further induction on~$\gamma_{1}$, one can prove
that
$a_{n} = \binom{n+\gamma_{1} - 1}{\gamma_{1} - 1}$.

Assume now that we have already proved the result for $\ell-1$: we want to prove it
for~$\ell$.
Denote $\tilde t := \tuple{t_{2}, \dotsc, t_{\ell}}$ and
%\Wlog, we may assume that $\gamma_{1} \geq 1$.
\[
%\begin{aligned}
% \gamma'_{i} &\coloneqq
% \begin{cases}
% \gamma_{i} & \text{if } i > 1\\
% \gamma_{1} - 1 & \text{if } i =1
% \end{cases}\\
g(\tilde t)  \coloneqq \frac{1}{\prod_{i=2}^{\ell}(1 - t_{i})^{\gamma_{i}}} 
\coloneqq \sum_{\tilde n \in \N^{\ell-1}} b_{\tilde n} \tilde t^{\, \tilde n}.
%\end{aligned}
\]
By inductive hypothesis, there exists $r(\tilde t) \in K[\tilde t]$ satisfying
(1) and (2) for~$g$.
Moreover, 
\[
f = g \cdot \frac{1}{(1 - t_{1})^{\gamma_{1}}} = 
\sum_{\tilde n \in \N^{\ell-1}} b_{\tilde n} \tilde t^{\, \tilde n} \cdot
\sum_{m \in \N} c_{m} t_{1}^{m},
\]
where $\frac 1 {(1 - t_{1})^{\gamma_{1}}} = \sum_{m \in \N} c_{m} t_{1}^{m}$.
Thus, denoting by $a_{m, \tilde n}$ the coefficient of 
$t_{1}^{m} \tilde t^{\, \tilde n}$ in~$f$, we have
\[
a_{m, \tilde n} = c_{m}b_{\tilde n}. 
\]
By the case $\ell = 1$ %and inductive hypothesis, we have that, for $\tilde n \in
%\N^{\ell -1}$ large enough, 
there exists $s(t_{1}) \in K[t_{1}]$ of degree $\gamma_{1}-1$
%(and actually equal to $\gamma_{1} -1$ if $K$ has characteristic~$0$)
such that, 
for every $m \in \N$ large enough,
\[
c_{m} = s(m).
\]
Thus, taking $m \in \N$ and $\tilde n \in \N^{\ell - 1}$ large enough, 
%and $m \geq m_{0}$,
we have
\[
a_{m,\tilde n} = s(m)r(\tilde n),
\]
and the polynomial $q(t_{1}, \tilde t) := s(t_{1})r(\tilde t) \in K[t]$
satisfies the conclusion.
% \[
% a_{m, \tilde n} = \Pa{\sum_{j=1}^{m_{0}-1} c_{j} 
% + \sum_{j=m_{0}}^{m} s(j)}  r(\tilde n)
% = \Pa{C + \sum_{j=m_{0}}^{m} s(j)}  r(\tilde n),
% \]
% where $C := \sum_{j=0}^{m_{0} -1} c_{j} \in K$.
% The second summand
% \[
% \sum_{j=m_{0}}^{m} s(j)  
% \]
% is equal to some polynomial $q'(\tilde n)$ for $m \geq m_{0}$, 
% such that, for
% For the first summamd, we have
% \[
% \sum_{j=0}^{m_{0}-1} c_{j} r(\tilde n) = C r(\tilde n),
% \]
% where $C := \sum_{j=0}^{m_{0} -1} c_{j} \in K$, and
% $C r(\tilde n)$ is also a polynomial in $(m, \tilde n)$ satisfying~(2).
\end{proof}

\begin{definition}\label{def:leader}
Let $p(\bar t) , q(\bar t)\in \R[\bar t]$; we write
$p = p_{0} + p_{1} + \dotsb + p_{d}$, where each $p_{i} \in \R[\bar t]$ is
homogeneous of degree $i$, and $p_{d} \ne 0$.  We call $p_{d}$ the 
\intro{leading homogeneous component} of~$p$ (if $p = 0$ then, by convention, the leading
homogeneous component of $p$ is $0$).  
As usual, if $\ell = 1$ we call the leading homogeneous component of $p$ the
leading \intro{term} of $p$.
\end{definition}
\begin{definition}\label{def:poly-order}
We write
\begin{itemize}
\item $p \preceq q$ if there exists $\cv \in \N^{\ell}$ such that, for every
$\bar n \in \N^{\ell}$ large enough, $p(\bar n) \le q(\bar n + \cv)$;
\item $p \simeq q$ if $p \preceq q $ and $q \preceq p$;
\item $p \succeq 0$ if, for every $\bar n \in \N^{\ell}$ large enough, $p(\bar n) \geq 0$.
\end{itemize}
\end{definition}
\begin{proposition}\label{prop:leader}
Let $p,q \in \R[\bar t]$ such that:
\begin{enumerate}
\item $p \succeq 0$,
\item $q \succeq 0$,
\item $p \simeq q$.
\end{enumerate}
Then, $p$ and $q$ have the same leading homogeneous component.
\end{proposition}
\begin{proof}
If either $p$ or $q$ is zero, it is clear that the other is also zero (and
therefore they have the same leading homogeneous component).
Thus, without loss of generality, we may assume that they are both non-zero.
Let $p'$ and $q'$ be the leading homogeneous component of $p$ and $q$ respectively,
and $h'$ be the leading homogeneous component of $h \coloneqq p - q$.
If, by contradiction,  $p' \ne q'$, then $\deg h = \max (\deg p, \deg q)$.
Let $\bar v \in \N^{\ell}$ such that $h'(\bar v) \ne 0$.
Then,
\[
r \coloneqq \lim_{s \to + \infty,\ s \in \N} h(s \bar v) \in \set{\pm \infty}.
\]
Since $q \preceq p$, we have $r = +\infty$, but since
since $p \preceq q$, we have $r = -\infty$, absurd.
\end{proof}

\section{Length functions and their entropy}\label{sec:length}

\subsection{Length functions}
% A \intro{length} function on the category $\Rmod$ of $R$-modules is a function
% \[
% \lambda: \Rmod \to \Rgeqzinf
% \]
% satisfying the following conditions:
% \begin{enumerate}
% \item $\lambda(0) = 0$;
% \item $\lambda(M) = \lambda(M')$ when $M$ and $M'$ are isomorphic;
% \item for every exact sequence
% $0 \to A \to B \to C \to 0$,
% \[
% \lambda(B) = \lambda(A) + \lambda(C);
% \]
% \item for every $M \in \Rmod$,
% \[
% \lambda(M) = \sup \set{\lambda(M'): M' \submod M\text{ finitely generated $R$-submodule}}
% \]
% \end{enumerate}

% Notice that we could  define a new function on $\Rmod$ by
% $\lambda'(N) := \sup\set{\lambda(M): M \submod N, \lambda(M) < \infty\rangle}$: then every $R$-module will be locally 
% $\lambda'$-finite, but unfortunately $\lambda'$ is not a length in general.

It can happen that $\lambda(R)$ is infinite: the following definition
deals with that situation.
\begin{definition}\label{def:lambda-finite}
Let $N$ be an $R$-module.
We say that $N$ is \intro{locally $\lR$-finite}%
\footnote{Also called ``locally $\lambda$-finite'' in \cite{SVV}: here we prefer to
write explicitly the ring $R$ too.}
if either of the following equivalent conditions hold:
\begin{enumerate}
\item For every $v \in N$, $\lambda(Rv) < \infty$;
\item For every $N_{0} \submod N$ finitely generated $R$-submodule, $\lambda(N_{0}) < \infty$.\end{enumerate}
\end{definition}

\begin{examples}\label{exs:length}
\begin{enumerate}[(a)]
\item If $R$ is a field, then the linear dimension is the unique length $\lambda$ on $\Rmod$
such that $\lambda(R) = 1$.
\item\label{ex:length-Z}  Let $R = \Z$, and define $\lambda(M)$
to be  the logarithm of the cardinality of $M$.
Then, $\lambda(R) = \infty$, and an Abelian group is
locally $\lR$-finite iff it is torsion.
We call $\lambda$ the \intro{standard length} on $\Z$-modules, and we will use it
often in examples.
\item Given any ring $R$,  the (classical) length of an $R$-module $M$ is the
 length of a composition series for $M$ (see \eg \cite{Eisenbud}).
\item\label{ex:trivial} The following are two ``trivial'' lengths:
\begin{enumerate}
\item $\lambda(M) = 0$ for every $M$;
\item  $\lambda(M) = \infty$ for every $M \neq 0$, and $\lambda(0) = 0$.
\end{enumerate}
\item\label{ex:singular} The following function is a length on $\Z$-modules:
\[
\lambda(A) \coloneqq
\begin{cases}
0 &\text{if $A$ is torsion}\\
\infty & \text{otherwise}.
\end{cases}
\]
\end{enumerate} 
\label{ex:length}
\end{examples}

%The following is an easy (and well known) observation, whose proof is left to
%the reader.
\begin{exercise}
$\lambda$ is nonzero  iff $\lambda(R) > 0$.
\end{exercise}

\begin{exercise}\label{exe:domain}
Let $R$ be an integral domain.
Then, there exists a unique length $\lambda_{0}$ on $R$-modules satisfying
$\lambda_{0}(R) = 1$.
Denoting by $K$ the field of fractions of $R$,
$\lambda_{0}$ 
 is defined by:
\[
\lambda_{0}(M) \coloneqq \dim_{K} (K \otimes M).
\]
\end{exercise}

The two trivial lengths in \ref{exs:length}(\ref{ex:trivial}) and the one in \ref{exs:length}(\ref{ex:singular}) are particular cases
of ``singular'' lengths, \ie lengths taking values only $0$ or $\infty$ (see
\cite[\S6]{Spirito} for a characterization): in the present treatment we will mostly ignore them, since the
associated entropies, the invariant $\mu$, and its modifications $\muhat$ and $\tildemu$ are all~$0$.

\medskip
\pagebreak[1]

A property we will use often in the rest of the paper is the following:
\begin{definition}\label{def:small}
Given an $R$-algebra  $T$ and a $T$-module $M$, we say that
$M$ is \mbox{\intro{$\lambda_T$-small}} if $M$ is finitely generated (as $T$-module)
and locally \lR-finite.%
\footnote{A similar notion is called ``Hilbert $T$-module'' in \cite[Ch.7]{Northcott:68}.}
\end{definition}
\begin{remark}
Assume that $T$ is a Noetherian $R$-algebra and let $M$ be a $T$-module.
The following are equivalent:
\begin{enumerate}
\item $M$ is $\lambda_{T}$-small;
\item there exists $V \submod M$ $R$-submodule such that:
\begin{enumerate}
\item $TV = M$,
\item $V$ is finitely generated (as $R$-module),
\item $\lambda(V) < \infty$;
\end{enumerate}
\item every submodule and every quotient of $M$ is $\lambda_{T}$-small.
\end{enumerate}
\end{remark}
Any submodule $V \submod M$ satisfying (2) in the above remark is a
\intro{witness} of the $\lambda_{T}$-smallness of~$M$.

% We conclude this subsection with a conjecture.
% \begin{conjecture}
% \begin{enumerate}
% \item Let $P \ideal R$ be a prime ideal.
% Then, there exists a unique length $\lambda_{P}$ on $\Rmod$ such that:
% \[\begin{aligned}
% \lambda_{P}(R/P) &= 1\\
% \lambda_{P}(M_{P}) &= \lambda(M)
% \end{aligned}\]
% where $M_{P}$ is the localization of $M$ at the prime~$P$.
% \item Let $\lambda$ be a length on $\Rmod$.
% Then, for every $P \ideal R$ prime ideal there exists 
% $r_{P}^{\lambda} \in R_{\geq 0} \cup \set{\infty}$ such that
% \[
% \lambda = \sum_{P \ideal R \text{ prime ideal}} r_{P}^{\lambda} \cdot \lambda_{P}.
% \]
% \item If $P \submod R$ is a maximal ideal, then $r_{P}^{\lambda}$ as in (2) is unique.
% \end{enumerate}
% \end{conjecture}

\cite[Thm.5]{Vamos:68} characterizes length functions on Noetherian rings: for
every prime ideal $P \submod R$ there is a canonical length function $l_{P}$ 
on~$\Rmod$, and any length function $\lambda$ can be written as
\[
\lambda = \sum_{P \ideal R \text{ prime ideal}} r_{P}^{\lambda} \cdot l_{P}
\]
for some $r_{P}^{\lambda} \in \R_{\geq 0} \cup \set \infty$ (we use the convention that
$0 \cdot \infty = 0$).

\subsection{Algebraic entropy}\label{sec:entropy}
The content of this and the following subsection can be skipped: it is mostly a motivation for
the definitions and results in the paper.
We recall the definition of algebraic entropy and its main properties.

%We fix  a length $\lambda$ on $\Rmod$.
Let $M$ be an $R$-module and $\phi$ be an endomorphism of $M$.
Given an $R$-submodule $V \submod M$, we define
\[
H_{\lambda}(\phi; V) := \lim_{n \to \infty} \frac{\lambda\Pa{V + \phi(V) + \dots + \phi^{n-1}(V)}}n
\] 
(the limit always exists by Fekete's Lemma, since the function
$n \mapsto \lambda\Pa{V + \phi(V) + \dots + \phi^{n-1}(V)}$ is subadditive: 
but see also later in this subsection).
The entropy of $\phi$ (according to the length $\lambda$) is defined by
\[
h_{\lambda}(\phi) = \sup\set{H_{\lambda}(\phi;V): V \submod M\ \ R\text{-submodule of finite length}}.
\]
Equivalently, we can see $M$ as an $R[X]$-module (with $X$ acting on $M$ as $\phi$),
and consider $h_{\lambda}$ as an invariant of $M$ as $R[X]$-module.
For the relationship between algebraic entropy and multiplicity, see
\cites{SVV, Northcott:68}.

More generally, given an $S$-module $M$, 
and an $R$-submodule $V \submod M$  of finite length,
%(where $S= R[x_{1}, \dotsc, x_{k}]$), let
%$S_{n}$ be the set of polynomials in $S$ of degree at most $n$.
define
\begin{align}\label{eq:H}
H_{\lambda}(M; V) & \coloneqq \lim_{n \to \infty} {\lambda(S_{n}V)}\big /{\tbinom {n+k}k}
= k! \lim_{n \to \infty} \lambda(S_{n}V) \big / n^{k} 
\\
h_{\lambda}(M) & \coloneqq \sup\set{H_{\lambda}(M;V): V \submod M\ \ R\text{-submodule of finite length}}.
\end{align}
The limit in the definition of $H_{\lambda}(M; V)$ exists, and $h_{\lambda}$ is the algebraic
entropy (relative to the length function $\lambda$).
We  prove the stronger result that $\lambda(S_{n}V)$ is eventually
equal to a polynomial as Theorem~\ref{thm:Hilbert-orbit}, 
and therefore the limit in \eqref{eq:H} exists.
However, the existence of the limit was already well-known:
\eg, \cites{CCK,DFG} give a more general
version. \cite{DFG}  consider 
the action of a cancellative amenable monoid:
in our case, we can identify $S$ with the group ring $R[\N^{k}]$, and therefore
an $S$-module is the same as an $R$-module $M $ together with an action $*$ of
$\N^{k}$ on $M$ by endomorphisms, and $\N^{k}$ is a cancellative amenable monoid.
Let $B_{n}$ be the set of tuples $\bar m \in \N^{k}$ such that
$m_{1} + \dotsb +  m_{k} \leq n$.
Then $\binom {n+k}k$ is the cardinality of  $B_{n}$. 
Moreover,  the family $(B_{n})_{n \in \N}$, is a F\o lner sequence for
$\N^{k}$
and therefore, as in \cite{DFG}, 
we can apply the machinery in \cite{CCK} to obtain that the
following limit exists (and is independent from the choice of the F\o lner sequence):
\[
\lim_{n \to \infty} \frac{\lambda(B_{n} * V) }{\card{B_{k}}} =
\lim_{n \to \infty} \frac{\lambda(S_{n} V) }{\tbinom {n+k}k}.
\]
See also~\cite{Virili:19} for an proof in the case when the acting monoid is a
finitely generated group.
% \footnote{Let $B^{k}_{n}$ be the set of monic monomials in $S_{n}$.
% Then $\binom {n+k}k$ is the cardinality of  $B^{k}_{n}$. 
% Moreover,  the family $(B^{k}_{n})_{n \in \N}$, is a F\o lner sequence
% inside the monoid of all monic monomials in $S$, and therefore we can apply
% the machinery of \cite{CCK}.}

\medskip
\pagebreak[1]

One of the most important properties of algebraic entropy is its additivity:
\begin{fact}\label{fact:add-entropy}
Let $0 \to A \to B \to C \to 0$ be an exact sequence of $S$-modules.
Assume that $B$ is locally $\lR$-finite.
Then,
\[
h_{\lambda}(B) = h_{\lambda}(A) + h_{\lambda}(C).
\]
\end{fact} 
%The above fact was already proved in \cite{SV:15, Virili:19, DFG}: we will show
%a stronger version in \ref{}.
We prove a stronger version of the above fact as Theorem~\ref{thm:additivity}.
However, the fact was well-known: see \cite{SV:15} for the case when $k=1$; 
\cite{DFG} gives a general version for $\Z$-modules 
 with the action of an amenable cancellative monoid (but the proof
 generalizes  to $R$-modules), while
\cite{Virili:19} treats the case of the action of an amenable finitely generated
group on $R$-modules.
\smallskip

Length functions on $R$-modules were explicitly introduced in \cite{NR:65} and further
studied in \cites{Vamos:PhD,  Vamos:68}; however, additive functions on
modules is a ``classical'' topic (see \eg\cite[Chapters 2 and 11]{AM}): one
of the novelties was allowing values in $\R_{\geq 0} \cup \set \infty$.
Algebraic entropy was introduced in \cite{DGSZ:09}: it and its variants have been
extensively studied, both in particular cases (\eg, $R = \Z$ and $k=1$) and in
general (including for non-Noetherian rings): see \eg\cites{SZ:09,
  DGZ:13, SV:15, BDGS1, GV:15, GS:17} and see \cites{SVV,GS:17b} for  surveys.

\subsection{Intrinsic algebraic entropy}\label{sec:entropy-intrinsic}
%We fix  a length $\lambda$ on $\Rmod$.
Let $M$ be an $R$-module and $\phi$ be an endomorphism of $M$.
Given an $R$-submodule $V \submod M$ such that
$\lambda((V + \phi(V)) / V) < \infty$, define
\[
\tilde H_{\lambda}(\phi; V) := \lim_{n \to \infty} \frac{\lambda\Pa{(V + \phi(V) + \dots + \phi^{n-1}(V)) / V}}n
\] 
(the limit always exists, again by Fekete's Lemma).
The intrinsic entropy of $\phi$ (according to the length $\lambda$), introduced in \cite{DGSV:15}, is defined by
\[
\tilde h_{\lambda}(\phi) = \sup\set{H_{\lambda}(\phi;V): V \submod M\ \ R\text{-submodule such that }
\lambda((V + \phi(V)) / V) < \infty
}.
\]
There is a corresponding addition theorem
\begin{fact}\label{fact:add-intrinsic-entropy}
Let $0 \to A \to B \to C \to 0$ be an exact sequence of $R[x]$-modules.
Then,
\[
\tilde h_{\lambda}(B) = \tilde h_{\lambda}(A) + \tilde h_{\lambda}(C).
\]
\end{fact} 
For a proof of the above fact, see \cites{DGSV:15, SV:18}: we will consider a
stronger version in a more general setting in \S\ref{sec:intrinsic}
(however, we were not able to prove additivity but only sub-additivity).

\section{Graded and filtered modules}\label{sec:filter}
%Let $S = R[x_{1}, \dotsc, x_{k}]$.

In this section we gather a few definitions and facts about graded and
filtered $S$-modules.
The most important ones are: how to construct a  graded
$S[y]$-module $\Rees(\overline A)$ starting from an upward filtered module $\overline A$ 
(Definitions~\ref{def:filtering} and \ref{def:blowup}),
and a  version of Artin-Rees Lemma for upward filtered
modules (Proposition~\ref{prop:AR-up}).

\subsection{Graded modules}
\begin{definition}%\ref{def:graded}
Fix $\bar \gamma = \tuple{\gamma_{1}, \dotsc, \gamma_{k}} \in \N^{k}$.
An \intro{$\N$-graded} $S$-module of degree $\bar \gamma$ is given by an $S$-module $M$ and a
decomposition 
\[
M = \bigoplus_{n \in \N} M_{n},
\]
where each $M_{n}$ is an $R$-module, and, for every $i \le k$ and $n \in \N$,
\[
x_{i} M_{n} \submod M_{n+ \gamma_{i}}.
\]
We denote by
$\overline M$ the module $M$ with the given grading (including the tuple  $\bar \gamma :=
\tuple{\gamma_{1}, \dotsc, \gamma_{k}}$).
\end{definition}

Given $\bj \in \N^{k}$ and $\bar \gamma \in  \N^{k}$, we denote
\[\begin{aligned}
\abs{\bj}_{\bar \gamma} & \coloneqq \bj \cdot \bar \gamma \coloneqq 
j_{1} \gamma_{1} + j_{2} \gamma_{2} + \dots + j_{k} \gamma_{k} \\
\x^{\bj} & \coloneqq x_{1}^{j_{1}} x_{2}^{j_{2}} \cdots x_{k}^{j_{k}}.
\end{aligned}
\]
Thus, if $\overline M$ is a graded module of degree $\bar \gamma$, then, for every
$\bj \in \N^{k}$ and $m \in M_{n}$,
\[
\x^{\bj} m \in M_{n + \abs{\bj}_{\bar \gamma}}
\]

We will use implicitly the following lemma many times in the remainder of the article.
\begin{lemma}\label{lem:fg}
Let $\overline M$ be an $\N$-graded $S$-module of degree $\bar \gamma$.
Assume that:
\begin{enumerate}
\item $\gamma_{\ell} > 0$ for each $\ell = 1, \dotsc, k$;
\item $M$ is finitely generated (as $S$-module).
\end{enumerate}
Then, each  $M_{n}$ is also finitely generated (as $R$-module).
\end{lemma}
\begin{proof}
Let $m_{1}, \dotsc, m_{p} \in M$ generate $M$ (as $S$-module).
Fix $n \in \N$; we want to show that  $M_{n}$ is finitely generated.
Without loss of generality, we may assume that each $m_{i}$ is homogeneous of degree $d_{i}$
(\ie, $m_{i} \in M_{d_{i}}$).

Let $a \in M_{n}$.
There exist $s_{1}, \dotsc, s_{p} \in S$ such that
\[
a = \sum_{i=1}^{p} s_{i} m_{i}.
\]
Write
\[
s_{i} = \sum_{\bj \in \N^{k}} r_{i, \bj}\, \x^{\bj} 
\]
for some (unique) $r_{i, \bj} \in R$.
Let 
\[
c_{i, \bj} \coloneqq x^{\bj} m_{i} \in 
S_{\abs{\bj}_{\bar \gamma}} M_{d_{i}} \subseteq M_{\abs{\bj}_{\bar \gamma} +  d_{i}}
\]
Thus,
\begin{equation}\label{eq:graded-decomposition}
a = \sum_{i=1}^{p} \sum_{\bj \in \N^{k}} r_{i,\bj}\, x^{\bj} m_{i}
= \sum_{i=1}^{p} \sum_{\bj \in \N^{k}} r_{i, \bj}\, c_{i,\bj}. 
\end{equation}
For every $\bj \in \N^{k}$, let
\begin{equation}\label{eq:wrong-sign}
%I_{\bj} \coloneqq \set{i \leq p:\ \abs{\bj}_{\bar \gamma} - d_{i} = n};
I_{\bj} \coloneqq \set{i \leq p:\ \abs{\bj}_{\bar \gamma} + d_{i} = n};
\end{equation}
notice that 
 $I_{\bj}$ is finite (since each $\gamma_{\ell} > 0$) 
and, for every $i \in I_{\bj}$, 
\[
c_{i,\bj} \in M_{n}.
\]

Since $a \in M_{n}$, we have that that in \eqref{eq:graded-decomposition} 
only the $c_{i, \bj}$ in $M_{n}$ 
%the summands with $r_{i, \bj}\, c_{i,\bj} \in M_{n}$ 
contribute to the sum:
that is, only the ones such that $i \in I_{\bj}$.
%and each 
%$c_{i, \bj} \in M_{\abs{\bj}_{\bar \gamma} + d_{i}}$,
Therefore,
\[
a = \sum_{\bj \in \N^{k}}\sum_{i \in I_{\bj}} r_{i,  \bj}\, c_{i,\bj} .
\]
Thus, $M_{n}$ is generated (as $R$-module) by the finite set
\[
\set{c_{i,\bj}: \bj \in \N^{k},\ \abs{\bj}_{\bar \gamma} \le n,
i \in I_{\bj}} \qedhere
\]
\end{proof}

\begin{definition}\label{def:acceptable-grading}
An $\N$-graded $S$-module $\overline A$ is \intro{acceptable} if:
\begin{enumerate}
\item $A$ is finitely generated (as $S$-module);
\item each $x_{i}$ has degree 1.
\end{enumerate}
\end{definition}

%For every $n \in \N$, we denote by $S^{(n)}$ the set of homogeneous polynomials
%in $S$ of degree $n$ (plus the $0$ polynomial).
%The following is an upward version of the Artin-Rees Lemma
%(however, notice that the proof works also when $R$ is non-commutative).
\begin{proposition}\label{prop:AR-graded-up}
Let $\overline A$ be an acceptable graded $S$-module.
Then, there exists $d \in \N$ such that, for every $n \in \N$,
\[
A_{d+n} = S^{(n)} A_{d}.
\]
\end{proposition}
\begin{proof}
It is always true that $S^{(n)} A_{d} \submod A_{d+n}$.
We want to show the opposite containment.

Let $a_{1}, \dotsc, a_{\ell} \in A$ be generators of~$A$.
Without loss of generality, 
we may assume that each $a_{i}$ is homogeneous of degree
$d_{i}$ (\ie, $a_{i} \in A_{d_{i}}$).

Let $d \coloneqq \max(d_{i}: i = 1, \dotsc, \ell)$.
Let $b \in A_{n+d}$.
We can write
\[
b = \sum_{i=1}^{\ell} s_{i} a_{i},
\]
for some $s_{i} \in S$.
For every $i = 1, \dotsc, \ell$, write
\[
s_{i} = \sum_{\bj \in \N^{k}} r_{i,\bj}\, \x^{\bj}
\]
for some (unique) $r_{i,\bj} \in R$.
Without loss of generality, as in the proof of Lemma~\ref{lem:fg}, 
we may assume that
$r_{i,\bj} = 0$ when $d_{i} + \abs{\bj} \ne n + d$.
For each $i, \bj$ such that $d_{i} + \abs{\bj} = n+d$, pick
$\bj[i], \bj''[i] \in \N^{k}$ such that:
\[
\bj'[i] + \bj''[i]= \bj \qquad \text{and} \qquad
\abs{\bj'[i]} =n. 
\]
Let 
\[\begin{aligned}
c_{i,\bj} &\coloneqq \x^{\bj'[i]} a_{i} \in S^{n-d_{i}} A_{d_{i}} \subseteq A_{n}\\
t_{i,\bj} &\coloneqq r_{i,\bj}\, \x^{\bj''[i]} \in S^{(n)}.
\end{aligned}\]
Thus,
\[
b = \sum_{i + \abs{\bj} =n} t_{i,\bj}\,c_{i,\bj} \in S^{(n)} A_{n}
\qedhere
\]
\end{proof}

\subsection{Filtered modules}
\begin{definition}\label{def:filtering}
Let $\bar \gamma := \tuple {\gamma_{1}, \dotsc, \gamma_{k}} \in \N^{k}$ and $N$ be an $S$-module.

An \intro{increasing filtering} on~$N$ with degrees~$\bar \gamma$ is an increasing sequence
of $R$\hyph{}submodules of $N$
\[
N_{0} \submod N_{1} \submod N_{2} \submod \dotsb \submod N
\]
%every $N_{i}$ is contained in $N_{i+1}$ and 
such that %$\bigcup_{i=0}^\infty N_{i} = N$, and 
$x_{i} N_{j} \submod N_{j + \gamma_{i}}$ for every $j \in \N$, $i \le k$.
We denote by $\overline N$ the $S$-module with the given tuple $\bar \gamma$ and the
filtering
$(N_{i})_{i \in \N}$.
Such a filtering is \intro{exhaustive} if $\bigcup_{i\in \N} N_{i} = N$.
\end{definition}

From now on, unless explicitly specified, all filterings will be \textbf{increasing}.

\begin{definition}\label{def:blowup}
% The ``Rees graded module'' associated to $\overline N$ is the following $S[y]$-module.
%\footnote{A similar object, the Rees algebra, is widely used in algebraic
%  geometry: see \eg \cite[\S6.5]{Eisenbud}.}
The \intro{blow-up} of the filtered $S$-module $\overline N$ is the following graded $S[y]$-module~$\ReesN$.%
\footnote{A similar construction is widely used in algebraic 
  geometry for \emph{downward} filtrations: see \eg \cite[\S5.2]{Eisenbud}.}

As a graded $R$-module, 
\[
\ReesN := \bigoplus_{n \in \N} N_{n} y^{n}.
\]
%The grading of $\ReesN$ is given by the decomposition 
%$\ReesN = \bigoplus_{n \in \N} N_{n} y^{n}$.

The multiplication by $x_{i}$ on $\ReesN$ is defined as:
\[
x_{i}(v y^{j}) := (x_{i}v) y^{j + \gamma_{i}},
\]
for every $i \le k$, $j \in \N$, $v \in N_{j}$, and then extended by $R$-linearity
on all $\ReesN$: notice that the  $x_{i}$   has degree $\gamma_{i}$ in~$\ReesN$.
The multiplication by $y$ on $\ReesN$ is defined as: 
\[
y (v y^{j}) := v y^{j+1},
\]
for every $j \in \N$, $v \in N_{j}$, and then extended by $R$-linearity
on all $\ReesN$: notice that $y$ has degree~$1$.
\end{definition}

\bigskip

%\subsection{Filtered modules}

Let $M$ be an $S$-module and
$\overline M = (M_{n})_{n\in \N}$ be a filtering of $M$ with degrees $\bar \gamma$.
For every $m \in \N$, we define 
\[
M^{m} := \bigoplus_{n \leq m} M_{n} y^{n} \le \Rees(M).
\]
We say that $M_{m}$ \intro{tightly generates} $\overline M$
%(as a filtered $S$-module) 
if:
for every $n \in \N$ and $v \in M_{n}$,
\begin{sentence}[($\dagger$)]
There exist $m_{1}, \dotsc, m_{r} \in \N$ with $m_{j} \leq m$, and
$v_{1}, \dotsc, v_{r} \in M$ such that $v_{j} \in M_{m_{j}}$,
and $\bar n_{1}, \dotsc, \bar n_{r} \in \N^{k}$ such that:
\[\begin{aligned}
v &= \bar x^{\bar n_{1}} v_{1} + \dots + \bar x^{\bar n_{r}} v_{r},\\
n &\ge \bar n_{j} \cdot \bar \gamma + m_{j}, \quad j = 1, \dots, r,
\end{aligned}
\]
\end{sentence}
where we
 are using the notations
\[
\bar n_{j} \cdot \bar \gamma := n_{j,1} \gamma_{1} + \dots + n_{j,k} \gamma_{k} \qquad \qquad
\x^{\bar n_{j}} := x_{1}^{n_{j,1}} \cdot \cdots \cdot x_{k}^{n_{j,k}}.
\]

Notice that ($\dagger$) is equivalent to:
\begin{sentence}[($\dagger'$)]
There exist $m_{1}, \dotsc, m_{r} \in \N$ with $m_{j} \leq m$,
$v_{1}, \dotsc, v_{r} \in M$ such that $v_{j} \in M_{m_{j}}$,
and  $p_{1}, \dotsc, p_{r} \in S$
such that: 
\[\begin{aligned}
v &= p_{1} v_{1} + \dots + p_{r}v_{r},\\
n &\ge \deg_{\bar \gamma}(p_{j}) + m_{j}, \quad j = 1, \dots, r .
\end{aligned}
\]
\end{sentence}

\begin{lemma}\label{lem:tight}
Let $m \in \N$.
$M^{m}$ generates $\Rees(\overline M)$ (as an $S[y]$-module) iff $M_{m}$ tightly
generates~$\overline M$. 
\end{lemma}
\begin{proof}
$\Rightarrow$)\ \
Let $n \in \N$ and $v \in M_{n}$.
%We want to find $m_{1}, \dotsc, m_{r} \in \N$,
%$v_{1}, \dotsc v_{r} \in M$,
%and $\bar n_{1}, \dotsc, \bar n_{r} \in \N^{k}$ 
%as in (\dagger).
Since $M^{m}$ generates $\RM$, there exist
$v_{1} y^{m_{1}}, \dotsc, v_{r}y^{m_{r}} \in M^{m}$ and
$q_{1}(\x,y), \dotsc, q_{r}(\x, y) \in R[\x,y]$ (remember that $S = R[\x]$)
such that:
\[
v y^{n} = q_{1}(\x, y) v_{1} y^{m_{1}} + \dots + q_{r}(\x, y) v_{r}y^{m_{r}}.
\]
Thus, if we define $p_{j}(\x) := q_{j}(\x, 0) \in S$, $j = 1, \dotsc, r$, we have
\[
v = p_{1} v_{1} + \dots + p_{r} v_{r}.
\]
Moreover,  $\deg_{\gamma}(p_{j}) + m_{j} \le m$,  $j = 1, \dotsc, r$, showing that
$M_{m}$ tightly generates~$\overline M$.

\smallskip

$\Leftarrow$)\ \
Let $n \in \N$ and $v y^{n} \in M_{n} y^{n}$.
Let $m_{1}, \dotsc, m_{r} \in \N$,
$v_{1}, \dotsc v_{r} \in M$ ,
and $\bar n_{1}, \dotsc, \bar n_{r} \in \N^{k}$ 
as in ($\dagger$).
For $j = 1, \dotsc, r$, define 
\[
d_{j} := n - (\bar n_{j} \cdot \bar \gamma + m_{j}) \in \N, \qquad
p_{j} := \x^{\bar n_{j}} y^{d_{j}} \in S[y]. 
\]
We have $v_{j} y^{m_{j}} \in M^{m}$ and
\[
v y^{n} = \sum_{i=1}^{m} p_{j } \cdot (v_{j} y^{d_{j}}) \in S[y] M^{m}.
\qedhere
\]
\end{proof}

\subsection{Acceptable filterings and upward Artin-Rees Lemma}
\begin{definition}\label{def:acceptable-filter}
Let $M$ be an $S$-module.
An \intro{acceptable filtering} of $M$ is given by a filtering $\overline M :=
(M_{n} : {n \in \N})$ such that $\Rees(\overline M)$ is an acceptable graded module;
that is:
\begin{enumerate}
\item each $x_{i}$ has degree 1;
%\item For every $n \in \N$, $\lambda(M_{n}) < \infty$;
\item $\Rees(\overline M)$ is finitely generated (as an $S[y]$-module).
%\item each $M_{i}$ is finitely generated (as $R$-module).
\setcounter{saveenum}{\value{enumi}}
\end{enumerate}
\end{definition}

The following is an upward version of Artin-Rees Lemma: however, as it can be
easily seen, the proof does not require $R$ to be a commutative ring.
\begin{proposition}\label{prop:AR-up}
Let $M$ be an $S$-module.
Let $\overline M$ be an exhaustive acceptable filtering of $M$.
Then, there exists $d \in \N$ such that:
\begin{enumerate}[(i)]
\item $M_{d}$ generates $M$ (as $S$-module);
\item for every $n \in \N$, $M_{n+d} = S_{n} M_{d}$.
\end{enumerate}
\end{proposition}
\begin{proof}
Since $\Rees(\overline M)$ is finitely generated, Lemma~\ref{lem:tight} implies
that there exists $d \in \N$ such that
$M_{d}$ tightly generates $\overline M$: thus, (i) is proven.

\begin{claim}
(ii) also holds (for the same $d$).
\end{claim}
By assumption, $\Rees(\overline M)$ is an acceptable graded $S[y]$-module.
Thus,
by (the proof of) Proposition~\ref{prop:AR-graded-up}, for every
$n \in \N$,
\[
M_{n+d} y^{n+d} = S[y]^{(n)} M_{d} y^{d}
\]
(as submodules of $\Rees(\overline M)$)
which is equivalent to (ii).
% The fact that $S_{n}M_{d} \subseteq S_{n+d}$ is clear by definition of filtering.
% We have to prove the opposite inclusion.
%
% Let $v \in M_{d+n}$.
% By ($\dagger$), there exist  exist $m_{1}, \dotsc, m_{r} \in \N$ with $m_{j} \leq d$,
% $v_{1}, \dotsc v_{r} \in M$ such that $v_{j} \in M_{m_{j}}$,
% and $\bar \ell_{1}, \dotsc, \bar \ell_{r} \in \N^{k}$ such that:
% \[\begin{aligned}
% v &= \bar x^{\bar \ell_{1}} v_{1} + \dots + \bar x^{\bar \ell_{r}} v_{r},\\
% d + n &\ge \abs{\bar \ell_{j}}  + m_{j}, \quad j = 1, \dots, r .
% \end{aligned}
% \]
% where we used the notation $\abs{\bar \ell} = \ell_{1} + \dots + \ell_{k}$.
% Split each $\ell_{j}$ into
% \[
% \ell_{j} = \ell_{j}' + \ell_{j''},
% \]
% in such a way that
% \[
% \abs{\ell_{j}'} \le n \ett \abs{\ell_{j}''} \le d - m_{j}.
% \]
% Define
% \[\begin{aligned}
% w_{j} &\coloneqq \x^{\bar \ell_{j}''} v_{j}  \in S_{d-m_{j}}M_{m_{j}} \subseteq M_{d}\\
% p_{j} &\coloneqq \x^{\bar \ell_{j}'} \in S_{\abs{\bar \ell_{j''}}} \subseteq S_{n}.
% \end{aligned}\]
% Thus, $v \in S_{n} M_{d}$.
\end{proof}

\section{Hilbert series for graded and filtered  modules}\label{sec:Hilbert}
%Let $S := R[x_{1}, \dotsc, x_{k}]$.
In this section we define the Hilbert series  associated to the length
function $\lambda$,  following the ideas in \cite{KLMP}
and~\cite[Ch.11]{AM}; in \S\ref{sec:Hilbert-p} we will define the
corresponding Hilbert polynomial.

\begin{thm}\label{thm:Hilbert-graded}
Let $\overline M$ be an $\N$-graded $S$-module of degree
$\bar \gamma \in \N^{k}$.
For every $n \in \N$, let $a_{n} := \lambda(M_{n})$.
Define
\[
F_{\overline M}(t) := \sum_{n} a_{n} t^{n}.
\]

Assume that:  
\begin{enumerate}
\item $\gamma_{i} > 0$ for $i = 1, \dotsc, k$;
\item $\lambda(M_{n}) < \infty$ for every $n \in \N$;
\item\label{as:Hilbert-graded:Noether}  $M$ is a Noetherian $S$-module.
\end{enumerate}
%and it is \lfgS.
%If $M$ is $\lambda$-Noetherian, 
Then, there exists a polynomial
$p(t) \in \R[t] $ such that
\[
F_{\overline M}(t) = \frac{p(t)}{\prod_{i=1}^{k}(1 - t^{\gamma_{i}})}.
\]
\end{thm}

\begin{proof}
By induction on $k$.

If $k = 0$, then, 
%since $M$ is $\lambda$-Noetherian, by Remark
%\ref{rem:lfgr-finite}, $\lambda(M)$ is finite, and therefore $\lambda(M_{n})$  is finite
%for every $n \in \N$.
%Moreover, 
since $M$ is Noetherian, only finitely many of the $M_{n}$ are
nonzero.
Thus, $F_{\overline M}(t)$ is a sum of finitely many (finite) terms, and hence it is
a polynomial.

Assume now that we have proven the conclusion for $k - 1$.
Let $y : M \to M$ be the multiplication by
$x_{k}$ and $\alpha := \gamma_{k}$.
For every $n \in \N$, let $y_{n} : M^{n} \to M^{n + \alpha}$ be the restriction of $y$
to $M_{n}$.
Let $K := \Ker(y)$ and $K_{n} := K \cap M_{n} = \Ker(y_{n})$.
Let $C_{n} := \Coker(y_{n}) = M_{m+ \alpha} / yM_{n}$, and
$C := \bigoplus_{n \in \N} C_{n}$.
Notice that both $\lambda(K_{n})$ and $\lambda(C_{n})$ are finite.
Therefore, both $K$ and $C$ are  $\N$-graded
$R[x_{1}, \dotsc, x_{k-1}]$-modules, and satisfy the assumptions of the
theorem
(that is, they are Noetherian modules, and each $K_{n}$ and each $C_{n}$ has
finite $\lambda$).

For every $n \in N$, consider the exact sequence
\[
0 \longto K_{n} \longto M_{n} \overset{y_{n}}\longto M_{n+ \alpha} \longto C_{n} \longto 0.  
\]
Since $\lambda$ is additive, we have
\[
a_{n + \alpha} - a_{n} = - \lambda(K_{n}) + \lambda(C_{n}). 
\]
Thus,
\[
\sum_{n} a_{n} t^{n + \alpha} - \sum_{n} a_{n} t^{n} = - \sum_{n} \lambda(K_{n})t^{n} + \sum_{n}  \lambda(C_{n})t^{n}.
\]

Therefore,
\[
(t^{\alpha}-1) F_{\overline M}(t) = - F_{\overline K} (t) + F_{\overline C}(t) 
\]
(where $\overline K$ is the $R[x_{1}, \dotsc, x_{k-1}]$-module with the given
grading, and similarly for $\overline C$).

Thus, by induction, there exist polynomials $q, q' \in \R[t]$ such that
\[
F_{\overline K}(t) = \frac{q(t)}{\prod_{i=1}^{k-1}(1 - t^{\gamma_{i}})}, \qquad
F_{\overline C}(t) = \frac{q'(t)}{\prod_{i=1}^{k-1}(1 - t^{\gamma_{i}})}.
\]
Therefore,
\[
F_{\overline M}(t) = \frac{q'(t) - q(t)}{\prod_{i=1}^{k}(1 - t^{\gamma_{i}})}.
\qedhere
\]
\end{proof}

\subsection{Filtered modules}

We move now from graded modules to (upward) filtered modules.

\begin{corollary}\label{cor:Hilbert-filtered}
Let $\overline N$
%$N = \bigcup N_{i}$ 
be a filtering on $N$ with degrees $\bar \gamma$.
%denoted by $\overline N$.

%For every $n \in \N$, let $a_{n} := \lambda(N_{n})$.
Define
\[
F_{\overline N}(t) := \sum_{n = 0}^{\infty} \lambda(N_{n}) t^{n}.
\]

Then,
\[
F_{\overline N}  = F_{\ReesN}.
\]

Therefore, if we assume that
\begin{enumerate}
\item $\gamma_{i} > 0$ for $i = 1, \dotsc, k$;
\item  $\lambda(N_{n}) < \infty$ for every $n \in \N$;
\item\label{ass:Hilbert-filtered:Noether} $\ReesN$ is Noetherian as $S[y]$-module;
\end{enumerate}
then, there exists a polynomial $p(t) \in \R[t]$ such that
\begin{equation}\label{eq:Hilbert-series-filter}
F_{\overline N}(t) = \frac{p(t)}{ (1 - t) \prod_{i = 1}^{k} (1 - t^{\gamma_{i}})}.
\end{equation}
\end{corollary}
\begin{proof}
Apply Theorem~\ref{thm:Hilbert-graded} to the graded ring
$\ReesN$.
The $(1 - t)$-factor in the denominator of \eqref{eq:Hilbert-series-filter}
is due to the action of $y$ on $\ReesN$ of degree~$1$.
\end{proof}

%{\texorpdfstring{$\lambda$}{\textlambda}-Hilbert  polynomials for Noetherian
%modules}
\section{Hilbert  polynomials for small modules}
\label{sec:Hilbert-p}
For the remainder of the article, excluding the appendices, we assume that $R$
is a \textbf{Noetherian} ring (commutative with $1$).
%Therefore, $S$ and $S[y]$ are also Noetherian rings.

\subsection{Hilbert polynomial for filtered modules}

\begin{definition}\label{def:good-filter}
Let $M$ be an $S$-module.
A \intro{good filtering} of $M$ is given by an acceptable filtering $\overline M :=
(M_{n} : {n \in \N})$ (see Definition~\ref{def:acceptable-filter})
such that:
% \begin{enumerate}
% \item Each $x_{i}$ has degree 1;
% %\item For every $n \in \N$, $\lambda(M_{n}) < \infty$;
% \item $\Rees(\overline M)$ is finitely generated (as an $S[y]$-module);
% \item each $M_{i}$ is finitely generated (as $R$-module).
% \setcounter{saveenum}{\value{enumi}}
% \end{enumerate}
% Such filtering is \intro{good} if moreover
\begin{enumerate}
\setcounter{enumi}{\value{saveenum}}
\item  $\forall n \in \N$ $\lambda(M_{n}) < \infty$.
%$\Rees(\overline M)$ is finitely generated (as an $S[y]$-module).
\end{enumerate}
\end{definition}

\begin{remark}
Let $\overline M$ be an acceptable filtering of an $S$-module $M$.
Then, since $R$ is Noetherian, $\Rees(\overline M)$ is Noetherian.
\end{remark}

\begin{thm}\label{thm:Hilbert-poly-filtered}
Let $M$ be an $S$-module and $\overline M := (M_{n})_{n \in \N}$ be a good filtering
of~$M$.
%For every $n \in \N$, let $a_{n} := \lambda(M_{n})$.
Then, for $n \in \N$ large enough, %$a_{n}$ 
$\lambda(M_{n})$ is equal to a polynomial 
$q_{\overline M}(n)$ of degree at most~$k$.

If moreover $\overline M$ is an exhaustive (and good) filtering of $M$, then the leading monomial of $q_{\overline M}$ does not depend on the choice of the
exhaustive good filtering (but only on $M$ and $\lambda$).
\end{thm}
Therefore, we can denote by $\mu(M)$ the leading monomial of the polynomial
$q_{\overline M}$ associated to some exhaustive good filtering of $M$ (if such good
filtering exists).
\begin{proof}
By Corollary~\ref{cor:Hilbert-filtered}, 
\[
F_{\overline M}(t) = \frac{p(t)}{(1-t)^{k+1}}, 
\]
for some polynomial $p(t) \in \R[t]$.
Thus, by Proposition~\ref{prop:poly-growth}, for $n$ large enough 
the coefficients $\lambda(M_{n})$ %$a_{n}$ 
of the power series $F_{\overline M}$ are equal to
some polynomial $q(n) \in \R[n]$ of degree at most~$k$.

Assume now $\overline M$ is exhaustive, 
and that $\overline M' := (M'_{n}: n \in \N)$ is another exhaustive good filtering of~$M$.
%and denote $a_{n}' := \lambda(M'_{n})$.
By Proposition~\ref{prop:AR-up}, there exists
$d_{0} \in \N$ such that, for every $n \in \N$,
\[
%M_{n + d_{0}} = S_{n} M_{d_{0}} \ett
M'_{n + d_{0}} = S_{n} M'_{d_{0}}
\]
% Since $\Rees(\overline M)$ is Noetherian, there exists $m_{0} \in \N$ such that
% $M^{m_{0}}$ generates $\Rees(M)$; similarly, after enlarging $m_{0}$ if
% necessary, we may assume that ${M'}^{m_{0}} := \bigoplus_{n \le m_{0}} M'_{n} y^{n}$
% generates $\Rees(M')$.
% Thus, $M_{m_{0}}$ tightly generates $\overline M$.
% Let $v_{1}, \dotsc, v_{\ell} \in {M'}^{m}$ which generate ${M'}^{m}$ (as
% $R$-module).
% Since $M = \bigcup_{n} M_{i}$, there exists $d \geq 0$ such that
% $v_{1}, \dotsc, v_{\ell} \in M_{m_{0} + d}$.
% Thus, 
% \[
% B'_{m_{0}} \leq B_{m_{0} + d}.
% \]
% \begin{claim}
% For every $n \geq m_{0}$,
% \[
% B'_{n} \submod B_{n + d}.
% \]
% \end{claim}
% Denote $S_{m}$ to be the set of polynomials in $S'$ of degree at most~$m$.
% Since $B'_{m_{0}}$ tightly generates $\overline M'$, for every $i \in \N$ we have
% \[
% B'_{m_{0} + i} \leq S_{i} B'_{m_{0}} \leq S_{i} B_{m_{0} + d} \le B'_{m_{0} + d + i}
% \]
% proving the claim.

Let $d_{1} \in \N$ such that $M'_{d_{0}} \le M_{d_{1}}$
($d_{1}$ exists because $M'_{d_{0}}$ is finitely generated as $R$-module).
Thus, for every $n$ large enough,
\[
q_{\overline M'}(n+d_{0}) = \lambda(M'_{n+d_{0}}) =
\lambda(S_{n}M'_{d_{0}}) \leq \lambda(S_{n}M_{d_{1}}) =
 q_{\overline M}(n + d_{1}).
\]
Similarly, 
\[
q_{\overline M}(n + d'_{0}) \le q_{\overline M'}(n + d'_{1})
\]
for some $d'_{0}, d'_{1} \in \N$ and every $n$ large enough, 
showing that $q_{\overline M}$ and $q_{\overline M'}$ have the same leading monomial.
\end{proof}

\subsection{Growth function}\label{sec:growth}

\begin{definition}
Let $N$ be an $S$-module, and $V_{0} \submod N$ be an $R$-submodule.
For every $n \in \N$, let 
%$S_{n}$ be the set of polynomial in $S$ of total degree at most $n$, 
$V_{n} := S_{n} V_{0}$ (notice 
that $V_{0} = S_{0} V_{0}$, and that $S_{n}$ and $V_n$ are $R$-modules).
We denote by
\[
\filt(V_{0}; N) \coloneqq (V_{n})_{n \in \N}
\]
the corresponding filtering of $N$ (as $S$-module), where each $x_{i}$ has
degree $1$,
and
\[
\gr(V_{0}; N) \coloneqq \Rees(\filt(V_{0}; N)) =
\bigoplus_{n \in \N} S_{n} V_{0}\, y^{n}
\]
be the corresponding graded  $S[y]$-module.
\end{definition}
Notice that $\gr(V_{0}; N)$ depends not on $N$ but only on $SV_{0} \submod N$.

\begin{thm}\label{thm:Hilbert-orbit}
Let $N$ be an $S$-module.
Let $V_{0} \submod N$ be an $R$-submodule.
% For every $n \in \N$, let $S_{n}$ be the set of polynomial in $S$ of total
% degree at most $n$, $V_{n} := S_{n} V_{0}$ (notice that $S_{0} =R$, 
% that $V_{0} = S_{0} V_{0}$, and that $S_{n}$ and $V_n$ are $R$-modules), and $a_{n} := \lambda(V_{n})$.
Define
\[
G_{V_{0}}(t) := \sum_{n \in \N} \lambda(S_{n} V_{0}) t^{n}. 
\]
Assume that:
\begin{enumerate}
\item $\lambda(V_{0}) < \infty$;
\item $V_{0}$ is finitely generated as $R$-module.
%\item $R$ is a Noetherian ring.
\end{enumerate}
Then, each $\lambda(S_{n} V_{0})$ is finite, and there exists a polynomial
$p(t) \in \R[t]$ such that
\[
G_{V_{0}}(t) = \frac{p(t)}{(1 - t)^{k+1}}.
\]
\end{thm}
\begin{proof}
First, we show that (1) implies that $\lambda(S_{n} V_{0})$ is finite for every $n \in \N$.
%where $V_{n} \coloneqq S_{n} V_{0}$.
In fact, $S_{n}V_{0}$ is a quotient of $V_{0}^{\ell}$, 
where $\ell \coloneqq \tbinom{n+k}n \in \N$ is the number of
monic monomials in $S$ of degree less or equal to $n$,
and $\lambda(V_{0}^{\ell})$ is finite.

%Let 
%\[
%V := SV_{0}.
%\]
Notice %that $V$ is an $S$-module, and 
that $\overline V \coloneqq \filt(V_{0};N)$
%$\overline V := \Pa{V_{n}: n \in \N}$ 
is an %exhaustive   
filtering of $N$ (as $S$-module).
%, where each $x_{i}$ has degree~$1$.
Moreover, $F_{\overline V}
%{\filt(V_{0};N)} 
= G_{V_{0}}$.
Thus, by Theorem~\ref{thm:Hilbert-poly-filtered}, it suffices to show that
$\Rees(\overline V) = \gr(V_{0}; N)$ is Noetherian as an $S[y]$-module to conclude
(since then the filtering $\overline V$ is good).
%Since, by (3),  $R$ is a Noetherian ring, $S[y]$ is also a Noetherian ring.
Since $S[y]$ is a Noetherian ring, it suffices to show that $\Rees(\overline V)$
 is finitely generated (as
an $S[y]$-module).
It is easy to see that $\Rees(\overline V)$  is generated by $V_{0} y^{0}$, and
the latter is finitely generated (as $R$-module) by (2).
\end{proof}

\begin{thm}[Hilbert polynomial]\label{thm:Hilbert-poly}
Let $N$ and $V_{0}$ %and $V_{n}$ 
be as in Theorem~\ref{thm:Hilbert-orbit} and 
assume that (1), (2) as in there hold.
Then, there exists a polynomial $q_{V_{0}}(t) \in \R[t]$ of degree at most~$k$, such that, for every $n$ large enough,
\[
\lambda(S_{n}V_{0}) = q_{V_{0}}(n).
\]

\textbf{Assume} moreover that
$N = SV_{0}$
(that is, $V_{0}$ witnesses that $N$ is  \lS-small).
Let $V_{0}'$ also  witness that $N$ is  \lS-small.
%finitely generated $R$-submodule of $N$ such that $\lambda(V_{0}') < \infty$ and
%$S V_{0}' = N$.
Then,  $q_{V_{0}}$ and $q_{V_{0}'}$ have the same leading term.

Therefore, if we define $\mu_{\lambda}(N)$ to be the leading term of $q_{V_{0}}$,
then $\mu_{\lambda}(N)$ does not depend on the choice of the witness~$V_{0}$.

Finally, let $c_{k}$ be the coefficient of $q_{V_{0}}$ of degree $k$. 
Then, 
\[
c_{k} = %k!\, h_{\lambda}(N) 
\frac{h_{\lambda}(N)}{k!} 
\] 
where $h_{\lambda}(N)$ is
the algebraic entropy of $N$ according to~$\lambda$ (see \S\ref{sec:entropy}).
\end{thm}
\begin{proof}
Let $\overline V \coloneqq (S_{n} V_{0})_{n \in \N}$.

Notice that $\overline V$ is a good filtering of~$N$, and that $\overline V$ is exhaustive
iff $SV_{0} = N$.
Therefore, Theorem~\ref{thm:Hilbert-poly-filtered} implies the existence of the
polynomial $q_{V_{0}}$, and that if $V_{0}$ witnesses that $N$ is $\lS$-small,
then the leading monomial of $q_{V_{0}}$ is independent of the choice of the
witness.

If $V_{0}$ is a witness, %let $c$ be the coefficient of degree $k$ of
%$q_{V_{0}}$.
then
\[
c_{k} = \lim_{n \to \infty} \frac{q(n)}{n^{k}} 
= \lim_{n \to \infty} \frac{\lambda(V_{n})}{n^{k}} 
%= \frac{H_{\lambda}(N,V_{0})}{k!} =
= \frac{h_{\lambda}(N)}{k!}. \qedhere
\]
%By Theorem~\ref{thm:Hilbert-orbit}.
% Given $V_{0}$ and $V_{0}'$ witnesses that $N$ is $\lS$-small,
% there exists a constant $n_{0}$ such that
% \[
% q_{V_{0}'}(t) \le q_{V_{0}}(t+n_{0}),
% \]
% and therefore, by and Proposition~\ref{prop:poly-growth}, the leading terms of
% $q_{V_{0}}$ and of $q_{V_{0}'}$ are the same.
\end{proof}

Notice that many authors (\eg, \cite{KLMP})
use a slightly different construction: in the
situation when $R$ is a field, they consider the function
\[
H(n) \coloneqq \dim(V_{n+1}/ V_{n}).
\]
It is easy to see that there exists a polynomial 
$\tilde G_{V_{0}}(t) \in \N[t]$ such that, for $n$ large enough,
$H(n) = \tilde G_{V_{0}}(t)$; from the definition it follows that
$\tilde G_{V_{0}}(t) = G_{V_{0}}(t+1) - G_{V_{0}}(t)$.
In the present situation,
we found it easier to work with the function~$\lambda(V_{n})$
(but see \S\ref{sec:intrinsic}).

\begin{definition}\label{def:mu-deg-size}
%Assume that $R$ is Noetherian, and 
Let $M$ be a  \lS-small $S$-module.
%such that $M$ is \lS-small (see Definition~\ref{def:small}).
%Thus, there exists $V_{0} \submod M$ witnessing that $M$
%$R$-submodule which is finitely generated, of
%finite length, and  generating $M$ as an $S$-module.
We define $\mu_{\lambda}(M)$ as in Theorem~\ref{thm:Hilbert-poly}
(with $\mu_{\lambda}(M) = 0$ iff $q_{V_{0}} = 0$): then,
$\mu_{\lambda}(M)$ does not depend on the choice of a witness. %$V_{0}$.
When $\lambda$ is clear from the context, we will write $\mu$ instead of $\mu_{\lambda}$.

Let $d$ be the degree of $\mu(M)$ and $m$ be the coefficient of $\mu(M)$.
We define the 
%\intro{degree} 
\intro{$\lambda$-dimension} of $M$ 
(as an $S$-module) to be equal to $d$,
 and its %\intro{$\lambda$-multiplicity}
\intro{$\lambda$-degree}
as 
\[
d!\, m %/ \tbinom {k+d} d
\]
When  $\lambda$ is clear, we will simply say ``dimension'' and
``degree'', respectively.%
\footnote{We boorrow the nomenclature from algebraic geometry, where, in the
  case when $R$ is a field and $\lambda$ is the linear dimension, the $\lambda$-dimension of
  $M$ is simply called the ``dimension'' of $M$, and the $\lambda$-degree is the
  ``degree'' of $M$: see \cite[\S1.9]{Eisenbud}. 
There is a different  construction in \cite[Ch.12]{Eisenbud} using
Hilbert-Samuel polynomial, where the ``degree'' becomes ``multiplicity'': see also \S\ref{sec:homogeneous}.}

If $\mu(M) = 0$, by convention we say that $M$ has dimension $-1$ and
degree undefined.
\end{definition}

One reason of the normalizing coefficient $d!$
%$1/\binom {k+d} d$ 
is the following:
\begin{example}
Assume that $0 < \lambda(R) < \infty$ and fix $d \leq k$.
Let  $M \coloneqq R[x_{1}, \dotsc, x_{d}]$ as an $S$-module, by
defining the action of $x_{j}$ on $M$ as multiplication by $0$ for $j > d$.
Let $V_{0} = R$ as a submodule of~$M$.
Then,
\[
q_{V_{0}}(n) = \binom{n+d}{d}\lambda(R) \qquad \mu(M) = \frac{t^{d}}{d!} \lambda(R)  
\]
Therefore, the $\lambda$-dimension of $M$ is $d$, and its  \ldeg
is $\lambda(R)$.
\end{example}

%\footnote{Maybe a better nomenklature would be
%  $\lambda$-dimension for the degree of $\mu$ and $\lambda$-multiplicity for the coefficient
% of~$\mu$} 

\begin{remark}\label{rem:dim-0}
Let $M$ be an $\lS$-small $S$-module.
Then, $\mu(M) = 0$ iff $\lambda(M) = 0$.
Moreover, the dimension of $M$ is $0$ iff $0 < \lambda(M) < \infty$.
\end{remark}

Notice that, if $\lambda(R) = 0$, then $\mu(M) = 0$ for every $\lS$-small $S$-module:
hence we often assume that $\lambda(R) \neq 0$ in the following.

\begin{proposition}\label{prop:degree}
Assume that 
$\lambda(R) = \ell$ with $0 < \ell < \infty$.
Let $p(\x) \in S$ of degree $e>0$.
Assume that $p_{e}$, the leading homogeneous component of $p$ (see Def.~\ref{def:leader}),
is not a zero-divisor (in~$S$).
Let $M \coloneqq S/(p)$.
Then, the \ldim of $M$ is $k-1$ and it \ldeg is $\ell e$.
\end{proposition}
\begin{proof}
``Usual'' proof.
Let $N := (p) \ideal S$.
We choose $M_{0} \coloneqq R \submod M$. %and $N_{0} \coloneqq (p) \cap S_{} \submod N$.
We have $S M_{0} = M$. %and $S N_{0} = N$.
For each $n \in \N$, we denote $M_{n} \coloneqq S_{n} M_{0} = S_{n} / ((p) \cap S_{n}) \submod M$.
%and $N_{n} \coloneqq S_{n} \cap (p) \submod N$.
Notice that, since $p_{e}$ is not a zero-divisor, for every $n \in \N$ we have
\[
p \cdot S_{n} = (p) \cap S_{n+e}
\] 
and therefore 
$S_{n+e} / p \cdot S_{n}$ and $M_{n+e}$
are isomorphic (as $R$-modules).
%we have $N_{n+e} = p \cdot N_{n}$, and
Moreover
also $p$ is not a zero-divisor, 
and therefore the multiplication by $p$ is an injective function (on~$S$).
%as $S$-modules.
Therefore, for every $n \in \N$,
the following sequence is exact:
\[
0 \longto S_{n} \overset{\cdot p}\longto S_{n+e} \longto M_{n+e} \longto 0
%0 \longto (p) \cap S_{n} \overset{\cdot p}\longto S_{n+e} 
\]
Therefore, if $q^{(S)}$ and $q^{(M)}$ are the Hilbert polynomials
associated to $S$ and $M$ respectively,
%to the filtrations $(N_{i})_{i \in \N}$ and $(M_{i})_{i \in \M}$, 
we have that, for every $n \in \N$ large enough,
\[
q^{(M)}(n+e) = q^{(S)}(n+e) - q^{(S)}(n).
\]
The conclusion follows.
\end{proof}
\begin{corollary}
Assume that  %$R$ is Noetherian and 
$\lambda(R) = 1$.
Let %$M$ and 
$p \in S$ be as in Proposition~\ref{prop:degree}.
Then, the \ldeg of $S/(p)$ is equal to $\deg(p)$, and in particular it is
independent from~$\lambda$. 
\end{corollary}
The above corollary implies that, if $R$ is an integral domain and $\lambda(R) = 1$,
then the \ldeg of $S/(p)$ does not depend on $\lambda$ (since the leading homogeneous
component of $p$ is not a zero divisor).
However, this is hardly surprising, since under the above assumption
$\lambda$ is unique (see Exercise~\ref{exe:domain}).

\subsection{Necessity of Noetherianity}\label{subsec:ex-noetherian}
We give an example of a $\lS$-small module over a non Noetherian ring $T$
with no associated Hilbert polynomial.

Define the following ring
\[
T \coloneqq \bigoplus_{p} \Z/p\Z 
\]
where $p$ varies among the set of primes.
Thus, $T$ is a direct sum of fields and it is not Noetherian (notice that it is
also not unitary).
Any $T$-module $M$ can be decomposed uniquely into the direct sum of its
$p$-components: 
\[
M = \bigoplus_{p} M_{p}
\]
where each $M_{p}$ is a $\Z/p\Z$-vector space with a certain dimension 
$\dim_{p}(M_{p})$.
Fix a sequence $\Pa{\alpha_{p}: p \text{ prime}}$ of real numbers such that, for each
prime $p$, $0 < \alpha_{p} < 1$, and $\sum_{p} \alpha_{p} = 1$.
Define the followng lenght function on $\Tmod$ as 
\[
\lambda(M) = \sum_{p} \alpha_{p} \dim_{p}(M_{p});
\]
thus, $\lambda(T) = 1$.
Let $M \coloneqq S \coloneqq T[x]$, where we see $S$ as a ring and $M$ as a $T$-module.
To give to $M$ a structure as $S$-module, we specify the action of $x$ on $M$ in
the following way:
\[
x \cdot (v_{p} x^{i}) \coloneqq
\begin{cases}
v_{p} x^{i+1} & \text{ if } p > i,\\
0 & \text{ otherwise,}
\end{cases}
\]
where $v_{p} \in \Z/p\Z$, and extend it by linearity to all~$M$.
Let $V_{0} \coloneqq T \submod M$.
Thus, $V_{0}$ is a finitely generated $T$-submodule of $M$ of finite length, but $\lambda(S_{n+1} V_{0}/S_{n} V_{0})$ is a strictly
decreasing sequence of real numbers in $(0,1)$, and therefore $\lambda(S_{n} V_{0})$ is not
eventually equal to any polynomial.

\smallskip
Notice that the algebraic entropy, \ie the limit 
$\lim_{n \to \infty} \lambda(S_{n} V_{0}) / n$, still exists.

\section{Dimension and degree: the general case}\label{sec:general}

We defined $\mu(M)$ when $M$ is a \lS-small $S$-module.
We will extend the definition to the case when $M$ is not necessarily 
$\lS$-small.
We need first to explain what is the range of $\mu$.

\subsection{The value monoid}\label{subsec:monoid}
(Remember that we fixed $k \in \N$).
We define the following ordered monoids $\vm$ and $\vmc$.
An element of $\vmc$ is either $0$ or a monomial
$r t^{d}$, where $r \in \R_{> 0} \cup \set \infty$ and $d \in \set{0, 1 ,\dotsc, k}$.
Given a monomial $0 \neq r t^{d} \in \vmc$, 
its degree is $d$ and its coefficient is $r$; for completeness we  define the
degree of $0$ to be $-1$.
$\vm$ is the subset of $\vmc$ given by $0$ and the monomials
with coefficient  which is not~$\infty$.

Remember that we follow the convention that $r + \infty = \infty$ for every
$r \in \R \cup \set \infty$.
The sum of two monomials in $\vmc$ is defined as
\[
r t^{n} \oplus s t^{m} =
\begin{cases}
r t^{n} & \text{if } n > m\\
s t^{m} & \text{if } n < m \\
(r + s) t^{n} & \text{if } n = m
\end{cases}
\]
and $0 + \mu = \mu$ for every  $\mu \in \vmc$.
We also define an ordering $\leq$ on $\vmc$ with the rule that 
\[
r t^{n} \leq s t^{m}
\]
iff either $n < m$ or $n = m$ and $r \leq s$, and $0 \leq \mu$ for every $\mu \in \vmc$.

With the above definitions, $\tuple{\vmc; \oplus, 0, \leq}$ is a 
commutative ordered monoid (with
$0$ the neutral element) and $\leq$ is a linear
ordering.
Moreover, $\vm$ is a submonoid of $\vmc$.

Notice that $\leq$ is a complete ordering on $\vmc$:
given $I \subseteq \vmc$, its supremum
$\sup(I)$ is  $0$ if $I$ is empty or $I = \set{0}$;
otherwise, $\sup(I)$
 is the monomial $r t^{d}$, where
\[\begin{aligned}
d &:= \max\set{d':  d' \text{ is the degree of some } \mu \in I} \in \set{0, \dotsc, k} \\
r &:= \sup \set{r' : r' t^{d} \in I} \in \R_{> 0} \cup \set \infty.
\end{aligned}\]
Moreover, $0$ is the minimum of $\vmc$ and $\infty t^{k}$ is its maximum, 
and $\infty t^{k}$ is an absorbing element: as an ordered set,
$\vmc$ is isomorphic to the real interval $[0,1]$.

\bigskip

We give now an equivalent description of the value monoid.
Let 
\[
P \coloneqq \set{p \in \R[t]: \deg(p) \leq k \text{ and } p(t) \geq 0 \text{ eventually}}.
\]
We endow $P$ with the (total) quasi-ordering  $\preceq$ defined in Def.~\ref{def:poly-order}, and  the
binary operation $+$ given by pointwise addition.
It is easy to see that $\tuple{P; +, 0, \preceq}$ is an ordered commutative monoid,
and that
the equivalence relation $\simeq$ on $P$ in Def.~\ref{def:poly-order} is compatible with the
structure of ordered monoid.
Therefore, $P/\mathord{\simeq}$ is also an ordered monoid (and
the induced quasi-ordering on $P/\mathord{\simeq}$ is a linear ordering).

\pagebreak[1]

Proposition~\ref{prop:leader} easily implies the following:%
\footnote{It is also quite easy to see it directly,
  since we are dealing with polynomials in 1 variable.} 
\begin{remark}
$P/\mathord{\simeq}$ is isomorphic to $\vm$ (as an ordered monoid), with the
canonical isomorphism given
by the function mapping the equivalence class of a polynomial $p$ to the leading
term of~$p$.
\end{remark}
\begin{remark}
$\vmc$ is the completion of $\vm$ (as an ordered set).
\end{remark}

\subsection{Non-small modules}\label{subsec:mu-general}
Let $M$ be an $S$-module (which might not be $\lS$-small).
Given $M' \submod M$ $S$-submodule which is \lS-small
%locally $\lR$-finite and finitely generated (as $S$-module), 
(see  Definition~\ref{def:small}),
let $\mu(M')$ be as in Definition~\ref{def:mu-deg-size}: 
notice that $\mu(M') \in \vm$. %(notice that $\mu(M') = 0$ iff $\lambda(M') = 0$).

Thus, we can \intro{define} $\mu(M) \in \vmc$ as the supremum of $\mu(M')$, where $M'$ varies among
all the possible $S$-submodules $M' \submod M$ which are
\lS-small.
% (notice that, by definition, $\mu(M) = 0$ if $0$ is the only $S$-submodule of
% $M$ of finite length).
We can then define as before the \ldim and \ldeg of $M$ as the coefficient (up
to a multiplicative constant) and the degree of $\mu(M)$, respectively:
the latter can be infinite.

From Remark~\ref{rem:dim-0} the following follows immediately.
\begin{remark}\label{rem:dim-0-general}
$\mu(M) = 0$ iff all  submodules of $M$ of finite length have length $0$.
In particular, if $0$ is the only submodule of $M$ of finite length, then $\mu(M)
= 0$.
\end{remark}

\medskip

An analogy  from geometry that might help the intuition is the
following. A  semi-algebraic set $X \subseteq \R^{k}$ has a  dimension 
$d \in \set{0, \dotsc, k}$ and a corresponding $d$-dimensional (Hausdorff) measure 
$r := \mathcal H_{d}(X) \in \R_{> 0} \cup \set \infty$: we could define $\mu(X) \coloneqq r t^{d} \in \vmc$ as the object
encapsulating both the dimension and the measure of~$X$ 
(with $\mu(X) = 0$ iff $X$ is empty).
The definition of $\oplus$ is such that if $X$ and $Y$ are disjoint manifolds,
then $\mu(X \cup Y) = \mu(X) \oplus \mu(Y)$; if $X$ and $Y$ are not necessarily disjoint, then
\begin{equation}\label{eq:manifolds}
\mu(X \cup Y) \oplus \mu(X \cap Y)= \mu(X) + \mu(Y).
\end{equation}

Thus, the $\lambda$-dimension of $M$ is the analogue of the dimension of $X$, and the 
$\lambda$-degree of $M$ is the analogue of the measure of~$X$.
We will see in \S\ref{sec:additivity} that \eqref{eq:manifolds}
%the area-coarea formula from geometric measure theory 
has an analogue for $\lS$-small modules: the additivity of~$\mu$.

\section{Addition theorem for exact sequences}\label{sec:additivity}
%Fix $k \in \N$ and let $S := R[x_{1}, \dotsc, x_{k}]$.
In this section we will prove the following Theorem.
\begin{thm}\label{thm:additivity}
Let $0 \longto  A \overset{\iota}\longto B \overset{\pi}\longto C \longto 0$ be an exact sequence of $S$-modules.
Assume that 
%\begin{enumerate}
%\item $R$ is Noetherian;
%\item 
$B$ is locally $\lR$-finite. %(see Definition~\ref{def:lambda-finite}).
%\end{enumerate}
% $A$, $B$, and $C$ 
% satisfy the assumptions of
% Theorem~\ref{thm:Hilbert-poly}.
Then, $\mu(B) = \mu(A) \oplus \mu(C)$ (see \ref{subsec:monoid} for the definition of $\oplus$).
\end{thm}
Notice that, under the  assumptions of the above theorem, also $A$ and $C$ are locally $\lR$-finite. %and finitely generated.
Notice moreover that $\mu(B)$ might have coefficient~$\infty$.

It is well-known that without the assumption that $B$ is locally $\lR$-finite,
the theorem may fail.
\begin{example}\label{ex:non-additivity}
Let $R = \Z$ with the standard length $\lambda$ (see
Example~\ref{ex:length}\ref{ex:length-Z}).
Let $A \coloneqq B \coloneqq \Z[x]$ and $C \coloneqq \Z/2\Z[X]$.
Let $\iota: A \to B$, $a \mapsto 2 a$ and let $\pi: B \to C$ be the canonical projection.
Then, 
$0 \longto  A \overset{\iota}\longto B \overset{\pi}\longto C \longto 0$ is an exact
sequence of $\Z[x]$-modules, but
$\mu(A) =\mu(B) = 0$ while $\mu(C) = \log 2 \cdot t$.
\end{example}

The main ingredient is the following proposition, which treats the case of
$\lS$-smallness (where in particular $\mu(B)$ has finite coefficient).
\begin{proposition}[Additivity]\label{prop:additivity}
Let $0 \longto  A \overset{\iota}\longto B \overset{\pi}\longto C \longto 0$ be an exact sequence of $S$-modules.
Assume that 
%\begin{enumerate}
%\item $R$ is Noetherian;
%\item 
$B$ is \lS-small. 
%locally $\lR$-finite 
%\item $B$ is finitely generated (as $S$-module).
%\end{enumerate}
% $A$, $B$, and $C$ 
% satisfy the assumptions of
% Theorem~\ref{thm:Hilbert-poly}.
Then, $\mu(B) = \mu(A) \oplus \mu(C)$.
\end{proposition}
Notice that, under the  assumptions of the above propositions, also $A$ and
$C$ are \lS-small.
%locally $\lR$-finite and finitely generated.

\begin{proof}%[Proof of Prop.~\ref{prop:additivity}]
Let $B_{0}$ be an $R$-submodule of $B$ such that $B_{0}$ is finitely
generated, $\lambda(B_{0}) < \infty$, and $S B_{0} = B$.
%Let  $\overline B = \Pa{B_{i}}_{i \in \N}$ be the filtering of $B$ associated to $B_{0}$,
%\ie $B_{i} := S_{i} B_{0}$.
For every $i \in N$, define 
\[\begin{aligned}
B_{i} &:= S_{i} B_{0}\\
A_{i} &:= \iota^{-1}(B_{i})\\
C_{i} &:= \pi(B_{i}),
\end{aligned}\]
and define $\overline B := (B_{i}: i \in \N)$, 
$\overline A := (A_{i}: i \in \N)$, 
and $\overline C := (C_{i}: i \in \N)$.
Notice that $\overline A$, $\overline B$, and $\overline C$ are good filterings of $A$, $B$,
and $C$, respectively (see Definiton~\ref{def:good-filter}).
Thus, by Theorem~\ref{thm:Hilbert-poly-filtered}, for $n \in \N$ large enough,
$\lambda(A_{n})=  q_{\overline A}(n)$ and $\mu(A)$ is the leading term of $q_{\overline A}$, 
and similarly for $B$ and~$C$.
Moreover, for every $n \in \N$, $\lambda(A_{n}) + \lambda(C_{n}) = \lambda(B_{n})$: thus,
\[
q_{\overline A} + q_{\overline C} = q_{\overline B},
\]
and therefore $\mu(A) \oplus \mu(C) = \mu(B)$.
\end{proof}

\begin{proof}[Proof of Thm.~\ref{thm:additivity}]
Since $B$ is locally $\lR$--finite, every submodule of $A$, $B$, or $C$ is
locally $\lR$-finite.

\begin{claim}
\[
\mu(B) \le \mu(A) \oplus \mu(C).
\]
\end{claim}
Let $B' \submod B$ be an $S$-submodule which is finitely generated.
Define 
\[
A' := \iota^{-1}(B'), \qquad C' := \pi(B').
\]
Notice that the sequence
\[
0 \longto A' \longto B' \longto C' \longto 0
\]
is exact, and therefore, by Proposition~\ref{prop:additivity},
\[
\mu(B') = \mu(A') \oplus \mu(C') \le \mu(A) \oplus \mu(C).
\]
Taking the supremum among all the $B'$, we get the Claim.

\begin{claim}
\[
\mu(B) \geq \mu(A) \oplus \mu(C).
\]
Let $A' \submod A$  and $C' \submod C$ be finitely generated 
$S$-submodules.
Since $C'$ is finitely generated and $\pi$ is surjective, there exists
$B' \le B$ finitely generated and such that $\pi(B') = C'$.
Define 
\[
B'' := B' + \iota(A'), \qquad A'' := \iota^{-1}(B'').
\]
We have that the sequence
\[
0 \longto A'' \longto B'' \longto C' \longto 0
\]
is exact, and $B''$ is finitely generated and locally $\lR$-finite.
Thus, by Proposition~\ref{prop:additivity},
\[
\mu(A') \oplus \mu(C') \le \mu(A'') \oplus \mu(C') = \mu(B'') \le \mu(B). 
\]
Taking the supremum on the left-hand side among all possible $A'$ and $C'$, we get
the Claim.\qedhere
\end{claim}
\end{proof}

\section{Modules over \texorpdfstring{$R$}{R}-algebrae}\label{sec:algebra}
Let $T$ be a finitely generated commutative $R$-algebra (therefore,
$T$ is Noetherian).
Let $M$ be a $T$-module.
We want to define the $\lambda$-dimension of $M$ as a $T$-module.

Fix $\gamma_{1}, \dotsc, \gamma_{k}$ generators of $T$ as $R$-algebra.
Equivalently, we fix a surjective homomorphism of $R$-algebrae
\[
\phi: S \to T 
%R[x_{1}, \dotsc, x_{k}] \to T
\]
and denote $\gamma_{i} \coloneqq \phi(x_{i})$, $i = 1, \dotsc, k$.
We can therefore see $M$ as an %$R[x_{1}, \dotsc, x_{k}]$
$S$-module, and we denote
it either by $\tuple{M; \phi}$ or by $\tuple{M; \bar \gamma}$.

We assume $M$ is  $\lT$-small.
%that $R$ is Noetherian, and 
%that there exists $M_{0} \submod M$
%$R$-submodule such that $TM_{0} = M$, $\lambda(M_{0})< \infty$, and $M_{0}$ 
%is finitely generated as $R$-module.

Thus, we can use the above data to compute $\mu(M; \phi)$ (which will depend on~$\phi$).

We prove now that, while the coefficient of $\mu$ may depend on $\phi$, its degree
does not. 
Thus, we can define the \ldim of $M$ (as a $T$-module)
as the degree of $\mu(M; \phi)$.

\begin{definition}
Spelling out all the assumptions, assume that:
\begin{enumerate}
\item $R$ is Noetherian;
\item $T$ is a finitely generated commutative $R$-algebra;
\item $M$ is a $T$-module;
\item there exists $M_{0} \submod M$ finitely generated $R$-submodule, such that
$\lambda(M_{0}) < \infty$ and $T M_{0} = M$.
\end{enumerate}
Then, we can define as before the \ldim of $M$ as a $T$-module, and
this dimension does not depend on the choice of $M_{0}$ or of $\phi$.
\end{definition}

\begin{examples}
Fix   some length function $\lambda$ on $R$ such that $\lambda(R) = 1$.
\begin{enumerate}[a)]
\item
Let
$M := T = R[z]$, 
$\bar \gamma := \tuple{z}$, $\bar \delta := \tuple{z,  z^{3}}$.
%$\bar \gamma := \tuple{1}$, $\bar \delta := \tuple{1,  3}$.
Thus, $\tuple{M; \bar \gamma}$ is $R[z]$ seen as $R[z]$-module with
the canonical action,
while $\tuple{M;\bar \delta}$ is $R[z]$ seen as $R[x_{1}, x_{2}]$-module,
with $x_{1}$ acting as multiplication by $z$ and $x_{2} $ as
multiplication by $z^{3}$.
Then, $\mu(M; \bar \gamma) = t$, while $\mu(M; \bar \delta) = 3 t$.
\item
Let  % $\Gamma := \Z$, 
$M := T = R[z, z^{-1}]$, 
$\bar \gamma := \tuple{z,z^{-1}}$, $\bar \delta := \tuple{z,  z^{-3}}$.
%$\bar \gamma := \tuple{1,-1}$, $\bar \delta := \tuple{1,  -3}$.
Thus, $\tuple{M; \bar \gamma}$ is $T$ seen as $R[x_{1},x_{2}]$-module 
with $x_{1}$ acting as multiplication by $z$, and $x_{2}$ as multiplication by $z^{-1}$,
while $\tuple{M;\bar \delta}$ is $T$ seen as $R[x_{1}, x_{2}]$-module,
with $x_{1}$ acting as multiplication by $z$ and $x_{2} $ as
multiplication by $z^{-3}$.
Then, $\mu(M; \bar \gamma) = 2t$, while $\mu(M; \bar \delta) = 4 t$.
\end{enumerate}
In both examples, 
we see that the two modules have different degrees, but have the same dimension.
\end{examples}

It remains to prove that the dimension of $\tuple{M; \phi}$ does not depend on the choice
of~$\phi$. 
It is clear that it suffices to prove the following:
\begin{thm}
Let $\bar \delta \in T^{k'}$ be another tuple of generators of~$T$. 
%and
%$\psi: R[x_{1}, \dotsc, x_{k'}] \to M$ be the corresponding surjective homomorphism
%of $R$-algebrae.
Then, $\tuple{M; \bar \gamma}$ and $\tuple{M; \bar \delta}$ have the same dimension.
\end{thm}
\begin{proof}
After exchanging the r\^oles of $\bar \gamma$ and $\bar \delta$ if necessary, we may assume that
$k \geq k'$.
After extending $\delta$ by setting $\delta_{i} = 0$ for $i \ge k'$, we may assume that
$k = k'$. 

%As usual, we denote $S := R[x_{1}, \dotsc, x_{k}]$ and by $S_{n}$
%the set of polynomials in $S$ of degree at most $n$.

We denote by $\psi: R[x_{1}, \dotsc, x_{k}] \to M$ the  surjective homomorphism
of $R$-algebrae corresponding to~$\bar \delta$ 
(and by $\phi$ the one corresponding to~$\bar \gamma$).

For every $n \in \N$, define
\[\begin{aligned}
T_{n}  &:= \phi(S_{n}),      & T_{n}'  &:= \psi(S_{n}),\\
M_{n} &:= T_{n}M_{0},  & M_{n}' &:= T_{n}' M_{0} .
\end{aligned}\]
Notice that both $(T_{n})_{n \in \N}$ and $(T_{n}')_{n \in \N}$ are %upward
filterings of $T$ as an $R$-module,  that
$(M_{n})_{n \in \N}$ and $(M_{n}')_{n \in \N}$ are %upward
filterings of $M$ as an $R$-module, and
that each $T_{n}$,  $T_{n}'$, $M_n$ and $M_{n}'$ are
finitely generated (as $R$-modules).

Moreover, 
$M_{0}$ %(resp., $M_{1}'$) 
generates both $\tuple{M; \bar \gamma}$ and 
$\tuple{M; \bar \gamma'}$ as $S$-modules.
%$(M_{n})_{n \in \N}$ (resp., $(M_{n}')_{n \in \N}$) is a good filterings
%of $\tuple{M; \bar \gamma}$ 
%(resp., $\tuple{M; \bar \gamma}$: see Definition~\ref{def:good-filter}): 
Thus, we can apply Theorem~\ref{thm:Hilbert-poly};
we denote by $q$ (resp., $q'$) the Hilbert polynomial of
$\tuple{M; \bar \gamma}$ (resp., of $\tuple{M; \bar \gamma'}$).

Let $c \in \N$ such that $\gamma_{1}, \dotsc, \gamma_{k} \in T_{c}'$.
Thus, $T_{1} \submod T_{c}'$.

The following is then clear
\begin{claim}
For every $n \in \N$, $T_{n} \submod (T_{c}')^{n}$
\end{claim}
Therefore, for every $n \in \N$
\[
M_{n} = T_{n} M_{0} \submod  (T_{c}')^{n} M_{0} \submod T_{cn}' M_{0} = M'_{cn}
\]
Therefore, for every $n \in \N$,
\[
\lambda(M_{n}) \le \lambda(M'_{cn}).
 \]
Therefore, for every $n$ large enough,
\[
q(n) = \lambda(M_{n}) \le \lambda(M'_{cn}) = q'(nc),
\]
proving that $\deg q \le \deg q'$.
Exchanging the r\^oles of $\phi$ and~$\phi'$, we see that $q$ and
$q'$ have the same degree.
\end{proof}

We end this section with a comparison between $\lambda$-dimension and Krull
dimension for affine rings.
\begin{lemma}\label{lem:Krull-affine}
Let $R$ be a field and $\lambda$ equal to the linear dimension (as $R$-vector
spaces).
Let $T$ be a finitely generated $R$-algebra.
Then, the $\lambda$-dimension and the Krull dimension of
$T$ coincide.
\end{lemma}
\begin{proof}
Let $d$ be equal to the Krull dimension of $T$.
By Noether Normalization (see \cite[Thm.13.3]{Eisenbud}, 
there exists an $R$-subalgebra $A \submod T$ such that:
\begin{enumerate}
\item $A$, as an $R$-algebra, is isomorphic to the polynomial ring
$R[y_{1}, \dotsc, y_{d}]$;
\item $T$ is finitely generated as $A$-module.
\end{enumerate}
Thus, $T$ and $A$ have the same Krull dimension $d$.
\begin{claim}
For any $T$ finitely generated $A$-module, $\lambdadim(T) \leq \lambdadim(A)$
(where  $\lambdadim$ denotes the $\lambda$-dimension). 
\end{claim}
In fact, $T$ is a quotient of $A^{n}$ (for some $n \in \N$), 
and Theorem~\ref{thm:additivity} implies that 
\[
\lambdadim(T) \leq \lambdadim(A^{n}) = \lambdadim(A).
\]

Since moreover $A \submod T$, we conclude that
$\lambdadim(T) = \lambdadim(A) = d$.
%Finally, $A$ has $\lambda$-dimension $d$.
\end{proof}
Lemma~\ref{lem:Krull-affine} answers positively and extends the conjecture in
\cite[Remark~5.9]{BDGS1}.
% where it was conjectured that
%  $0 < h^{(1)}(T) < \infty$  if and only if the Krull dimension of $T$ is 1 
% (the notation $h^{(1)}$ will be introduced in \S\ref{sec:receptive}).

Here is another example of equality between Krull and $\lambda$ dimensions.
%(besides Lemma~\ref{lem:Krull-affine}).

\begin{lemma}\label{lem:Krull-discrete}
Fix $2 \le n \in \N$.
Let $R \coloneqq \Z/(n)$.
Let $\lambda$ be the length on $\Rmod$ given by
\[
\lambda(M) \coloneqq \log(\card M).
\]
Let $A$ be a finitely generated $R$-algebra.
%$A := S/J$ for some ideal $J  \ideal S$.
Then, $\lambdadim(A)$ is equal to the Krull dimension of $A$.
\end{lemma}
\begin{proof}
We denote by $\dim_{K}$ the Krull dimension. 
%and by $\lambdadim$ the $\lambda$-dimension.

Write the factorization of $n$ into primes:
\[
n = p_{1}^{e_{1}} \cdots p_{\ell}^{e_{\ell}}.
\]
Decompose $A$ into a direct sum of $Z_{p_{i}^{e_{i}}}$-algebrae~$A_{i}$.
Since $\dim(A) = \max(\dim(A_{i}: i = 1, \dotsc, \ell)$, where $\dim$ is either 
$\dim_{K}$ or $\lambdadim$,
%the Krull dimension or the $\lambda$-dimension, 
it suffices to treat the case when $\ell = 1$, \ie
$n = p^{e}$.

Let $B_{i} := p^{i} A$ as $R$-submodule of $A$, for $i = 1, \dotsc, e$. 
We have $0 = B_{e} \submod B_{e-1} \submod \dotsc \submod B_{0} = A$.
Notice that, for every $i < e$, $\tilde B_{i} \coloneqq B_{i}/B_{i+1}$ is a
$Z/(p)$-algebra.
Moreover, there exists a surjective homomorphism of $\Z$-modules between $\tilde
B_{0}$ and $\tilde B_{i}$, mapping $a + pA$ to $p^{i} a + p^{i+1} A$.
Thus, $\mu(\tilde B_{i}) \leq \mu(\tilde B_{0})$ and therefore 
\[
\lambdadim(A) = \lambdadim(\tilde B_{0}) = \lambdadim(A /pA).
\]
Moreover, for $\Z/(p)$-modules, $\lambda$ (up to a constant factor)
is equal to the linear
dimension, and Lemma~\ref{lem:Krull-affine} implies that
\[
\lambdadim(A/pA) = \dim_{K}(A/pA).
\]
Finally, $pA$ is the unique minimal prime ideal of $A$, and therefore
\[
\dim_{K}(A) = \dim_{K}(A/pA).
\qedhere
\]
% We proceed by induction on $n$.
% If $n$ is prime, then $R$ is a field, and $\lambda$ (up to a constant factor)
% is equal to the linear
% dimension:
% thus, the conclusion follows from Lemma~\ref{lem:Krull-affine}.

% Otherwise, let $p \in \N$ be a prime divisor of $n$ and let $m \coloneqq n/p$.
% We have the exact sequence of $Z/(n)$-modules
% \[
% 0 \longto p A \longto A \longto A/pA \longto 0.
% \]
% Both $pA$ and $A/pA$ are quotients of $\Z/(m)[\x]$: therefore, by inductive
% hypothesis, their $\lambda$-dimensions are equal to their Krull dimensions.
% Since $\dim(A) = \max(\dim(A/pA), \dim(pA))$ when $\dim$ is either the Krull
% dimension or the $\lambda$-dimension,
% the conclusion follows. 
\end{proof}

\section{Hilbert-Samuel polynomial for homogeneous modules}\label{sec:homogeneous}

Let %$S \coloneqq R[x_{1}, \dotsc, x_{n}]$ and 
$I \coloneqq (x_{1}, \dotsc, x_{k}) \ideal S$.
Let $M$ be an $S$-module.
For every $n \in \N$,  define $c_{n} \coloneqq \lambda(M/I^{n+1}M)$.
\begin{thm}
Assume that:
\begin{enumerate}[(i)]
%\item $R$ is Noetherian;
\item $M$ is finitely generated (as $S$-module);
\item $\lambda(M/IM)$ is finite.
\end{enumerate}
Then, for every $n \in \N$, $c_{n}$ is finite, and there exists a polynomial
$\bar q(t) \in \R[t]$ such that:
\begin{enumerate}
\item for every $n \in \N$ large enough,
$c_{n} = \bar q(n)$;
\item $\deg \bar q \le k$.
\end{enumerate}
\end{thm}
\begin{proof}
Usual proof (see e.g.\ \cite[Prop.12.2]{Eisenbud}).
\end{proof}
Assume moreover, besides the hypothesis in the theorem, that $V \submod M$
witnesses that $M$ is  \lS-small.
% there exists
% $V \submod M$ $R$-submodule, such that:
% \begin{enumerate}[(a)]
% \item $\lambda(V)$ is finite;
% \item $V$ is finitely generated (as $R$-module);
% \item $SV = M$.
% \end{enumerate}
%Then, one can define $a_{n} := \lambda(S_{n }V)$. %as before.
Notice that
\[
c_{n} = \lambda(M/I^{n+1}M) = \lambda(S_{n}V/I^{n+1}M) \le \lambda(S_{n }V). %a_{n}.
\]
Therefore, denoting by $q_{V}$ the Hilbert polynomial associated to $V$,
we have $\bar q(t) \leq q_{V}(t)$ for every $t$ large enough.
If $\bar \mu(M)$ is the leading term of $\bar q$, we have therefore $\bar \mu(M) \le
\mu(M)$.

In general, it can happen that $\bar \mu(M) < \mu(M)$.
\begin{example}
Let $K$ be a field, $\lambda$ be the linear dimension over $K$,
$S \coloneqq K[x_{1},x_{2}]$, $M \coloneqq K[x_{1},x_{2}]/(x_{1}x_{2} - 1)$.
Then, $\mu(M) = 2t$, while $\bar \mu(M) = 0$.
\end{example}

It is easy to prove that for homogeneous ideals the situation is different.
\begin{exercise}
%Assume that $R$ is Noetherian.
Let $J \ideal S$ be a homogeneous ideal, and $M \coloneqq S/J$.
Then,
\[
\mu(M) = \bar \mu(M).
\]
More precisely, fix a finite  set $G$ generating  $J$, and
let $n_{0}$ be the maximum degree of the polynomials in~$G$.
Let $V \coloneqq R$.
Then, for every $n > n_{0}$, $S_{n} V$ and $M/I^{n+1}$ are isomorphic (as
$R$-modules), and therefore 
\[
%a_{n} = 
\lambda(S_{n} V) = \lambda(M/I^{n+1}) = c_{n}.
\]
\end{exercise}
% \begin{question}
% Assume that:
% \begin{enumerate}
% \item $\lambda(R) < \infty$
% \item $R$ is Noetherian
% \item $J \submod S$ is a homogeneous ideal.
% \end{enumerate}
% Is it true that $\mu(S/J) = \bar \mu(S/J)$?
% \end{question}

See also \cite[Ch.12]{Eisenbud} and \cite[Ch.7]{Northcott:68}
for the ``classical'' version of the Hilbert-Samuel polynomial.

%A similar notion for the algebraic entropy of a single endomorphism
%was defined in \cite{DGSV:15}.

\section{\texorpdfstring{$d$}{d}-dimensional and receptive versions of entropy}\label{sec:receptive}
%Let $S,M,V_{0}$ be as in Definition~\ref{def:mu-deg-size}.
Let %$R$ be Noetherian and 
$M$ be an 
%\lS-small  ?
$S$-module.
Let $m$ be the coefficient of $\mu(M)$.
For every $d \le k$, define
\[
h^{(d)}_{\lambda}(M) \coloneqq
\begin{cases}
\infty & \text{if } \deg \mu(M) > d\\
0  & \text{if } \deg \mu(M) < d\\
%\frac m {d!} 
d!\, m & \text{if } \deg \mu(M) = d.
\end{cases}
\]

% Assume now that $M$ satisfies the following weaker condition:
% \begin{enumerate}
% \item $R$ is Noetherian.
% \end{enumerate}
% Define:
% \[
% h^{(d)}_{\lambda}(M) \coloneqq \sup \set{h^{(d)}(M'):
% M' \submod M \text{ \lS-small $S$-submodule}}
% \]
The value $h^{(1)}(M)$ is the receptive entropy of $M$ \wrt 
the standard regular system generated by
$(x_{1}, \dotsc, x_{k})$ (see \cites{BDGS1, BDGS2});
we call each $h^{(d)} (M)$ the \intro{$d$-dimensional   entropy} of~$M$
(and thus the algebraic entropy $h$ is the $k$-dimensional  entropy).

The case $d= 1$ of the following Proposition answers positively  (and extends) \cite[Question~5.10]{BDGS1}.
\begin{proposition}\label{prop:intrinsic-entropy-additive}
Let $0 \to A \to B \to C \to 0$ be an exact sequence of $S$-modules.
Assume that 
%\begin{enumerate}
%\item  $R$ is Noetherian;
%\item 
$B$ is locally \lR-finite.
%\end{enumerate}
Then, for every $d \le k$ %\deg \mu(B),
\[
h^{(d)} (B) = h^{(d)}(A) + h^{(d)}(C). 
\]
\end{proposition}
\begin{proof}
Denote by $\dim(A)$ the $\lambda$-dimension of $A$ (that is, the degree of $\mu_{\lambda}(A)$).
By Theorem~\ref{thm:additivity}, $\dim(B) = \max(\dim(A), \dim(C))$.
Assume, for simplicity, that $\dim(A) \leq \dim(C)$, and therefore
$\dim(B) = \dim(C)$
 (the other case when $\dim(A) > \dim(C)$ is similar).

If $d < \dim(B)$, then $h^{(d)}(B) = h^{d}(C) = h^{(d)}(A) = 0$.

If $d < \dim(B)$, then $h^{(d)}(B) = h^{d}(C) = \infty$.

If $d = \dim(B) = \dim(A)$, then  Theorem~\ref{thm:additivity} again implies
that $h^{(d)} (B) = h^{(d)}(A) + h^{(d)}(C)$.

If $d = \dim(B) > \dim(A)$, then $h^{(d)}(A) = 0$ and 
Theorem~\ref{thm:additivity} again implies that
$h^{(d)}(B) = h^{(d)}(C)$.

In all four cases, the conclusion follows. 
\end{proof}

\medskip
Let $T$ be a finitely generated $R$-algebra (thus, $T$ is Noetherian).
We can give similar definitions of entropies for $T$-modules.
Fix $\bar \gamma = \tuple{\gamma_{1}, \dotsc, \gamma_{k}}$ generators of $T$ (as  $R$-algebra).
Let $\tuple {M, \bar \gamma}$ be the $S$-algebra defined in \S\ref{sec:algebra}.
%Assume that $M$ satisfies the conditions in \S\ref{sec:algebra}.
% Therefore, $\tuple {M, \bar \gamma}$ satisfies the conditions in
% Definition~\ref{def:mu-deg-size}, and we can therefore define, for every $d \le
% \deg \mu(M;\bar \gamma)$,
% \[
% h^{(d, \bar \gamma)}_{\lambda}(M) \coloneqq h^{(d)}_{\lambda}(M; \bar \gamma).
% \]
%For general $M$ (under the assumption that $R$ is Noetherian) we define
%Assume that $R$ is Noetherian. 
We define
\[
h^{(d, \bar \gamma)}_{\lambda}(M) \coloneqq \sup \set{
h^{(d)}(M'; \bar \gamma):
M' \submod M %\text{ locally \lR-finite, finitely generated $T$-submodule}
\text{ \lT-small $T$-submodule}
}.
\]

$h^{(1, \bar \gamma)}_{\lambda}(M)$ is the receptive entropy of $M$ \wrt 
the standard regular system generated by
$\bar \gamma$ (see \cites{BDGS1, BDGS2});
we call each $h^{(d, \bar \gamma)} (M)$ the \intro{$d$-dimensional   entropy} of~$M$
\wrt $\bar \gamma$.

\begin{thm}
Let $0 \to A \to B \to C \to 0$ be an exact sequence of $T$-modules.
Assume that 
%\begin{enumerate}
%\item  $R$ is Noetherian;
%\item 
$B$ is locally \lR-finite.
%\end{enumerate}
Then, for every $d \le k$, % \deg \mu(B),
\[
h^{(d, \bar \gamma)} (B) = h^{(d, \bar \gamma)}(A) + h^{(d, \bar \gamma)}(C). 
\]
\end{thm}
In particular, the receptive  entropy $h^{(1, \bar \gamma)}$
is additive (under the assumptions of
Noetherianity of $R$ and local \lR-finiteness!).

\section{Totally additive versions of \texorpdfstring{$\mu$}{\textmu} and (receptive) entropy}\label{sec:hat}
The definition of $\mu(M)$  reflects the usual definition of algebraic
entropy (see \S\ref{sec:entropy}).
Following a  construction in \cite[Prop.3]{Vamos:68},
we propose an alternative invariant, which is in some ways better behaved.

\begin{definition}
Let $A$ be an $S$-module.
A $\lS$-small chain in $A$ is a sequence of $S$-submodules
\[
\Aseq = \Pa{A_{1} \submod A_{2} \submod \dots \submod A_{2n-1} \submod A_{2n} \submod A},
\]
where, for every $i \le n$, 
\[
\hat A_{i} \coloneqq A_{2i}/A_{2i-1}
\]
is \lS-small.
%locally $\lR$-finite and finitely generated.
We call $n$ is the \intro{size} of $\Aseq$.
\end{definition}
\begin{definition}
Let $\theta$ be a partial function from $\Smod$
to~$\vmc$.
%$S$-modules
% to the commutative monoid
% of monomials of the form $m t^{d}$, where $d \in \N$ and
% $m \in R_{\geq 0} \cup \set {\infty}$ (examples of such function are  
% $\lambda$  and~$\mu$).

We will be interested only in functions $\theta$ which satisfy the following
conditions:
\begin{description}
\item[Domain] the domain of $\theta$ includes all $\lS$-small $S$-modules;
\item[Additivity] 
$\theta(0) = 0$ and,
for every exact sequence
$0 \to A \to B \to C \to 0$ of 
$\lS$-small
$S$-modules, % in the domain of $\theta$,
$\theta(B) = \theta(A) \oplus \theta(C)$;
\item[Invariance]
if $A$ and $B$ are isomorphic $S$-modules in the domain of $\theta$,
$\theta(A) = \theta(B)$.
\end{description}
%Assume that the domain of $\theta$ includes all the \lS-small $S$-modules.
Let $A$ be any $S$-module.
Given a $\lS$-small chain $\Aseq$ in $A$ of size $n$, we define
\[\begin{aligned}
\theta(\Aseq) &\coloneqq \sum_{i=1}^{n} \theta(\hat A_{i})\\
\thetahat(A) &\coloneqq \sup\set{\theta(\Aseq): \text{ $\Aseq$ $\lS$-small chain in $A$}}.
\end{aligned}\]
\end{definition}
% We will see in \S\ref{sec:additivity-muhat} that $\muhat$ satisfies a stronger
% additivity condition than~$\mu$.

We will see later that $\thetahat$ can be defined in a simpler way
(Proposition~\ref{prop:cosmall}); see also \cite[Prop.3]{Vamos:68} for an equivalent approach.

% \begin{question}
% In the definition of $\thetahat$, can we restrict ourselves to consider only
% chains on size 1?
% \end{question}

% \begin{remark}
% $\muhat \geq \mu$.
% \end{remark}

For ``well-behaved'' length functions $\lambda$, we have $\hat \lambda = \lambda$ (here we take
$S = R$).
However, the following example shows that it is not always the case.
\begin{example}
Let $\lambda$ be any singular non-zero length (\eg, the length in Example~\ref{exs:length}(\ref{ex:singular})).
% following function on $\Z$-modules:
% \[
% \lambda(A) \coloneqq
% \begin{cases}
% 0 &\text{if $A$ is torsion}\\
% \infty & \text{otherwise}.
% \end{cases}
% \]
Then, $\lambda \neq 0$ but $\hat \lambda = 0$.
\end{example}

\begin{proposition}\label{prop:theta}
Assume:
\begin{enumerate}
%\item $R$ is Noetherian;
\item the domain of $\theta$ includes all $\lS$-small $S$-modules;
\item $\theta$ is additive and invariant (on $\lS$-small $S$-modules). 
\end{enumerate}
Then, 
\begin{enumerate}[(a)]
\item $\thetahat$ is also additive and invariant;
\item if $A$ is a \lS-small $S$-module, then $\thetahat(A) = \theta(A)$;
\item 
For every $S$-module $A$,
\[
\thetahat(A) = \sup \set{\thetahat(B): B \submod A \text{ finitely generated $S$-submodule}}
\]
\item if $A$ is a locally \lR-finite $S$-module, then
\[
\thetahat(A) = \sup \set{\theta(B): B \submod A \text{ finitely generated $S$-submodule}}
\]
\item
\[
\hat{\thetahat} = \thetahat.
\]
\end{enumerate}
\end{proposition}
\begin{proof}
The proof is quite straightforward; we will prove that $\thetahat$ is
additive, and leave the remainder as an exercise (see also \cite{Vamos:68}).
Thus, let $0 \longto  A \overset{\iota}\longto B \overset{\pi}\longto C \longto 0$ be an exact sequence of $S$-modules.
\begin{claim}
\[
\thetahat(B) \le \thetahat(A) \oplus \thetahat(C).
\]
\end{claim}
Let $\Bseq = \Pa{B_{1} \submod B_{2} \submod \dots \submod B_{2n-1} \submod B_{2n} \submod B}$
be a $\lS$-small chain in $B$.
For every $i \le n$, define
\[
A_{i} \coloneqq \iota^{-1}(B_{i}) \qquad C_{i} \coloneqq \pi(A_{i}).
\] 
Then, 
$\Aseq \coloneqq \Pa{A_{1} \submod A_{2} \submod \dots \submod A_{2n-1} \submod A_{2n} \submod A}$
and $\Cseq \coloneqq \Pa{C_{1} \submod C_{2} \submod \dots \submod C_{2m-1} \submod C_{2m} \submod C}$
are $\lS$-small chains in $A$ and $C$, respectively.
Moreover, for every $i \leq n$, we have an exact sequence
\[
0 \longto \hat A_{i} \longto \hat B_{i} \longto \hat C_{i} \longto 0.
\]
Therefore, $\theta(\Bseq) = \theta(\Aseq) \oplus \theta(\Cseq)$, and the claim follows.

\begin{claim}
\[
\thetahat(B) \geq \thetahat(A) \oplus \thetahat(C)
\]
\end{claim}
Let $\Aseq = \Pa{A_{1} \submod A_{2} \submod \dots \submod A_{2n-1} \submod A_{2n} \submod A}$
and $\Cseq = \Pa{C_{1} \submod C_{2} \submod \dots \submod C_{2m-1} \submod C_{2m} \submod C}$
be $\lS$-small chains in $A$ and $C$, respectively.
For every $i \le 2(m+n)$
\[
B_{i} \coloneqq
\begin{cases}
\iota(A_{i}) & \text{if } i \le 2n;\\
\pi^{-1}(C_{i-2n}) &\text{if } 2n < i \leq 2(n+m).
\end{cases}
\]
Then, $\Bseq \coloneqq \Pa{B_{1} \submod B_{2} \submod \dots \submod B_{2n+2m-1} \submod B_{2n + 2m} \submod B}$
is a $\lS$-small chain in~$B$.
Moreover, for every $i \leq n+m$,
\[
\hat B_{i} =
\begin{cases}
\hat A_{i} & \text{if } i \leq n;\\
\hat C_{i} & \text{if } n < i \leq n+m.
\end{cases}
\]
Therefore, $\theta(\Bseq) = \theta(\Aseq) \oplus \theta(\Cseq)$, and the claim follows.
\end{proof}

%For the remainder of this section, we \textbf{assume} that $R$ is \textbf{Noetherian}.

\begin{definition}
Given an ideal $I \ideal R$, we say that $I$ is \lcosmall if
$\lambda(R/I) < \infty$.
\end{definition}
\begin{remark}
Let $I \ideal R$ be a \lcosmall ideal, and $A$ be an $R$-module.
Then, $A/IA$ is locally \lR-finite.
\end{remark}

\begin{proposition}\label{prop:cosmall}
Let $\theta$ be as above and total.
Assume that, for every $S$-module~$A$
\[
\theta(A) = \sup \set{\theta(B): B \submod A \text{ \lS-small $S$-submodule}}.
\]
Then,
\[
%\begin{aligned}
\thetahat(A) \geq \theta(A)
\]
and, if $A$ is finitely generated,
\[
\thetahat(A) = \sup \set{\theta(A/IA): I \ideal R \text{ \lcosmall ideal}}.
%\end{aligned}
\]
\end{proposition}
The proof of the above proposition is in the next subsection: for now we will
record some consequences.

% \begin{corollary}
% $\lambda = \hat \lambda$ iff
% \[
% \lambda(R) = \sup \set{\lambda(A/I): I \ideal R }
% \]
% \end{corollary}

\begin{corollary}\label{cor:muhat}
\begin{enumerate}
\item $\muhat$ satisfies the conclusions of Propositions~\ref{prop:theta}
and~\ref{prop:cosmall};
\item $\muhat(S) = \hat\lambda(R) t^{k}$;
\item $\muhat(A) = \mu(A)$ for every locally \lR-finite $S$-module~$A$;
\item if $\lambda(R) < \infty$, then $\muhat = \mu$.
\end{enumerate}
\end{corollary}

\begin{corollary}
For every $d \le k$,
\begin{enumerate}
\item  the $d$-dimensional entropy (see \S\ref{sec:receptive})
$\hat{h}^{(d)}_{\lambda}$ is a length functions (on all $S$-modules) and satisfies
the conclusion of Proposition~\ref{prop:cosmall};
\item $\hat{h}^{(d)}_{\lambda}(A) = h^{(d)}_{\lambda}(A)$ for every locally \lR-finite $S$-module~$A$;
\item if $\lambda(R) < \infty$, then $\hat{h}^{(d)}_{\lambda} = h^{(d)}_{\lambda}$.
\end{enumerate}
\end{corollary}

Remember that
\[
\hat h_{\lambda} = \hat h^{(k)}
\]
and therefore from the above Corollary we obtain that $\hat h$ is a length
function $\Smod$, that $\hat h(A) = h(A)$ when $A$ is locally finite, and
$\hat h = h$ when $\lambda(R) < \infty$.

\begin{corollary}
Let $T$ be a finitely generated $R$-algebra, 
$\bar \gamma \in T^{k}$ be a set of generators of $T$.
Given $d \le k$, let $h^{(d, \bar \gamma)}$
be defined as in \S\ref{sec:receptive}.
Then:
\begin{enumerate}
\item 
$\hat{h}^{(d, \bar \gamma)}_{\lambda}$ is a length functions (on all $T$-modules) and satisfies
the conclusion of Proposition~\ref{prop:cosmall};
\item $\hat{h}^{(d, \bar \gamma)}_{\lambda}(A) = h^{(d, \bar \gamma)}_{\lambda}(A)$ for every locally \lR-finite $S$-module~$A$;
\item if $\lambda(R) < \infty$, then $\hat{h}^{(d, \bar \gamma)}_{\lambda} = h^{(d, \bar \gamma)}_{\lambda}$.
\end{enumerate}
\end{corollary}
\begin{proof}
Apply Proposition~\ref{prop:theta} to the function $h^{(d, \bar \gamma)}$.
\end{proof}

\begin{corollary}
Let $\lambda$ be the standard length on $\Z$-modules introduced in 
Example~\ref{ex:length}(\ref{ex:length-Z}).
Then, for every finitely generated $\Z[\x]$-module~$A$,
\[
\muhat(A) = \sup \set{\mu(A/n A): 2 \leq n \in \N} =
\lim_{n \to \infty } \mu(A/n!\, A).
\]
\end{corollary}

We cannot drop the assumption that $A$ is finitely generated in
Proposition~\ref{prop:cosmall}.
\begin{example}\label{ex:hat-Q}
Let $R \coloneqq \Z$ and $\lambda$ be the standard length.
% introduced in Example~\ref{ex:length}\ref{ex:length-Z}.
Let $A \coloneqq \Q[\x]$ (seen as a $\Z[\x]$-module).
Then, 
\[\begin{aligned}
\muhat(A) &= \infty \cdot t^{k} \\
\sup \set{\mu(A/nA): 2 \leq n \in \N} &= 0.
\end{aligned}\]
\end{example}

\subsection{Proof of Proposition~\ref{prop:cosmall}}

\begin{lemma}\label{lem:cosmall-interesection}
Let $I, J \ideal  R$ be \lcosmall ideals.
Then, $I \cap J$ is also \lcosmall.
\end{lemma}
\begin{proof}
$R/I\cap J$ embeds into  $R/I \times R/J$.
\end{proof}

\begin{lemma}
Let $I, J \ideal R$ be \lcosmall ideals.
Then, $IJ$ is also \lcosmall.
\end{lemma}
\begin{proof}
Let $\av = (a_{1}, \dotsc, a_{\ell})$ generate~$I$.
Then, $I/IJ$ is a quotient of $(R/J)^{\ell}$ via the map
\[
(r_{1} + J, \dotsc, r_{\ell} + J) \mapsto r_{1} a_{1} + \dots + r_{\ell} a_{\ell} + IJ.
\]
Therefore, $\lambda(I/IJ) \le \ell \lambda(R/J) < \infty$,
and
\[
\lambda(R/IJ) \le \lambda(R/I) + \lambda(I/IJ) < \infty. \qedhere
\]
\end{proof}

\begin{lemma}
Let $A$ be an $S$-module.
Assume that $A$ is $\lS$-small.
Then, 
\[
Ann_{R}(A) \coloneqq \set{r \in R: r A = 0}
\]
 is \lcosmall
\end{lemma}
\begin{proof}
Let $a_{1} ,\dotsc, a_{\ell}$ be generators of $A$ (as $S$-module).
For every $i \le \ell$, $\lambda(Ra_{i}) < \infty$.
Moreover, $Ra_{i}$ is isomorphic (as $R$-module) to
$R/Ann_{R}(a_{i})$, and therefore $Ann_{R}(a_{i})$ is \lcosmall.
Finally,
\[
Ann_{R}(A) = Ann_{R}(a_{1}) \cap \cdots \cap Ann_{R}(a_{\ell})  
\]
and the conclusion follows from  Lemma~\ref{lem:cosmall-interesection}.
\end{proof}

\begin{proof}[Proof of Proposition~\ref{prop:cosmall}]
Let $B \submod A$ be an \lS-small $S$-submodule.
By definition, $\thetahat(A) \geq \theta (B)$: therefore, 
by the assumption, $\thetahat(A) \geq \theta(A)$.

\smallskip

Define 
\[
\theta'(A) \coloneqq \sup\set{\theta(A/IA): I \ideal R\ \text{\lcosmall ideal}}.
\]
We want to prove that, when $A$ is finitely generated, $\thetahat(A) = \theta'(A)$.
It suffices to show that $\theta'$ is additive on finitely generated $S$-modules.
Thus, let
\[
0 \longto A \longto B \longto C \longto 0
\]
be an exact sequence of finitely generated $S$-modules.
\begin{claim}
$\thetahat(B) \le \thetahat(A) \oplus \thetahat(C)$.
\end{claim}
Let $I \ideal R$ be a \lcosmall ideal.
We have the exact sequence of \lS-small $S$-modules
\[
0 \longto A/(A \cap IB) \longto B/IB \longto C/IC \longto 0.
\]
Since $\theta$ is additive on \lS-small $S$-modules, and $IA \submod A \cap IB$,
we have
\[
\thetahat(A) + \thetahat(C) 
\geq \theta(A/IA) \oplus \theta(C/IC) 
\geq \theta(A/(A \cap IB)) \oplus \theta(C/IC)
= \theta(B/IB),
\]
and the claim follows.

\begin{claim}
$\thetahat(B) \geq \thetahat(A) \oplus \thetahat(C)$.
\end{claim}
Let $I, I' \ideal R$ be \lcosmall ideals.
We want to prove that
\[
\theta(A/IA) + \theta(C/I'C) \leq \thetahat(B).
\]
Replacing $I, I'$ with $I \cap I'$, 
without loss of generality we may assume that $I = I'$.
By Artin-Rees Lemma (see \eg \cite[\S4.7]{Northcott:68}), there exists
$1 \leq n_{0} \in \N$ such that, for every $m \in N$,
\begin{equation}
A \cap I^{m + n_{0}} B = I^{m}(A \cap I^{n_{0}} B).
\label{eq:Artin-Rees}
\end{equation}
Let $J \coloneqq I^{n_{0}}$ and $A' \coloneqq A \cap JB$: notice that $J$ is also \lcosmall.
Taking $m \coloneqq n_{0}$ in~\eqref{eq:Artin-Rees}, we obtain:
\[
A \cap J^{2} B = J A'.
\]
Thus, we have the exact sequence
\[
0 \longto A/ JA' \longto B/J^{2}B \longto C/J^{2}C \longto 0.
\]
The modules appearing above are all \lS-small: therefore,
\[
\theta(A/IA) \oplus \theta(C/IC) \le \theta(A/JA') \oplus \theta(C/J^{2}C) = \theta(B/J^{2}B) \le \thetahat(B),
\]
proving the Claim.
%
% Let $\Aseq = (A_{1} \submod A_{2} \submod \dots \submod A_{2n})$ be a \lS-small chain in~$A$.
% For each $i = 1, \dotsc, n$, let $J_{i} \coloneqq Ann_{R}(A_{2i}/A_{2i-1})$.
% By the above lemma, each $J_{i}$ is a \lcosmall ideal, and therefore
% $J \coloneqq J_{1} \cap \dots \cap J_{n}$ is \lcosmall.
% We have
% \[
% \theta(\Aseq) \le \theta(A/J)
% \]
% proving that
% \[
% \thetahat(A) \le \sup \set{\theta(A/I): I \ideal R \text{ \lcosmall ideal}}.
% \]
%
% \begin{claim}
% If $B$ is a locally \lR-finite $S$-module, then $\theta(B) = \thetahat(B)$.
% \end{claim}
% In fact,
% \[
% \theta(B) = \sup\set{\theta(C): C \leq B \text{ \lS-small $S$-submodule}}
% \le \thetahat(B).
% \]
%
% Moreover,
% \[
% \thetahat(A) \geq \sup \set{\thetahat(A/I): I \ideal R \text{ \lcosmall ideal}}
% \]
% and, for every $I \ideal R$ \lcosmall ideal,
% $\thetahat(A/I) \geq \theta(A/I)$.
\end{proof}

% \begin{example}
% Let $R \coloneqq \Z$ and $\lambda$ be as in
% Example~\ref{ex:length}(\ref{ex:length-Z}).
% Let $M \coloneqq S \coloneqq Z[x]$.
% Then, $\mu(M) = 0$, while $\muhat(M) = \infty\cdot t^{1}$.
% \end{example}

\subsection{Examples}
Let $R \coloneqq \Z$, $\alpha$ be the standard length introduced in
Example~\ref{ex:length}(\ref{ex:length-Z}) and $\beta$ be the length given by the
rank (\ie, $\beta(M) = \dim_{\Q}(M \otimes \Q)$).
Since $\beta(\Z) = 1 < \infty$, $\mu_{\beta} = \muhat_{\beta}$.

1)
Let  $S \coloneqq \Z[x]$,  $I$ be an ideal of $S$, and $M := S/I$.
The following table shows the values of $\mu_{\alpha}(M)$,  $\muhat_{\alpha}(M)$,
and $\mu_{\beta}(M)$ for some values of~$I$:

\medskip

\begin{tabular}{r@{}llll}
\toprule
$I$      &                & $\mu_{\alpha}(S/I)$    & $\muhat_{\alpha}(S/I)$ & $\mu_{\beta}(S/I)$  \\
\midrule
$0$      &                & $0$             & $\infty \cdot t^1$         & $t^1$         \\
$S$      &                & $0$             & $0$               & $0$           \\
$(n)$    &; $2 \le n \in \N$ & $\log(n) \cdot t^1$ & $\log(n) \cdot t^1$   & $0$           \\
$(p(x))$ &; $\deg p \geq 1$ & $0$             & $\infty \cdot t^{0}$       & $\deg p\cdot t^0$ \\
$(p(x), n)$ &; $2 \le n \in \N$, &  \multirow{2}{*}{$\log(n) \deg p \cdot t^0$} & 
   \multirow{2}{*}{$\log(n) \deg p \cdot t^0$ } & \multirow{2}{*}{$0$} \\
  \multicolumn{2}{r}{$p$~monic, $\deg p \geq 1\,$} \\
\bottomrule
\end{tabular}

\bigskip

2) Let $S \coloneqq \Z[x_{1}, x_{2}]$ and  $M \coloneqq S/(x_{1}x_{2})$.
Then, 
\[
\mu_{\alpha}(M) = 0, \quad \muhat(M) = \infty \cdot t, \quad \mu_{\beta}(M) = 2 t.
\]

\section{Intrinsic Hilbert polynomial}\label{sec:intrinsic}

%In this section, we assume that $R$ is Noetherian.

In \cite{DGSV:15} the authors introduced the ``intrinsic'' algebraic entropy, a
variant of the more usual algebraic entropy: following a similar pattern, we
introduce here the intrinsic Hilbert polynomial.

Let $A$ be an $S$-module.
Let $\overline A \coloneqq (A_{i})_{i \in \N}$ be a  filtering on~$A$.
% (see Definition~\ref{def:good-filter}).
For each $i \in \N$, define $\tilde A_{i} \coloneqq A_{i+1}/A_{i}$ (as
$R$-modules).
Define
\[
\btA \coloneqq \bigoplus_{i \in \N} \tilde A_{i} 
\]
as graded $S$-module (where all the $x_{i}$ have degree 1).
An equivalent description of $\btA$ is the following.
Remember that $\Rees(\overline A)$ is an $S[y]$-module.
Let 
\[
\cdot y: \Rees(\overline A) \to \Rees(\overline A)
\]
 be the multiplication by $y$.
Then,
\[
\btA = \Coker(\cdot y);
\]
notice that $\Coker(\cdot y)$ is an $S[y]$-module: however, $y$ acts trivially on
$\Coker(\cdot y)$, hence we lose nothing in considering $\Coker(\cdot y)$ as an
$S$-module; moreover, the above isomorphism is of graded $S$-modules.
In particular, if $\Rees(\overline A)$ is Noetherian (as $S[y]$-module), then
$\btA$ is also Noetherian (as $S$-module).

\begin{definition}
We say that $\overline A$ is a \intro{$\lambda$-inert filtering} on $A$ if:
\begin{enumerate}
\item $\btA$ is Noetherian (as $S$-module);
%\item $\Rees(\overline A)$ is Noetherian (as $S[y]$-module);
%\item each $A_{i}$ is finitely generated (as $R$-module); 
\item there exists
$n_{0} \in \N$ such that, for every $n \geq n_{0}$, 
\[
\lambda(\tilde A_{n}) < \infty.
\]
\end{enumerate}
\end{definition}

\begin{proposition}\label{prop:intrinsic-Hilbert}
Assume that 
%$\Rees(\overline A)$ is Noetherian (as $S[y]$-module) and that there exists
%$n_{0} \in \N$ such that, for every $n \geq n_{0}$, 
%\[
%\tilde a_{n} \coloneqq \lambda(\tilde A_{n})
%\]
%is finite.
$\overline A$ is a $\lambda$-inert filtering on $A$.
Then, there exists a polynomial $\tilde q_{\overline A}(t) \in \Q[t]$ of degree at
most $k-1$ such that, for every $n \in \N$ large enough,
\[
\lambda(\tilde A_{n}) = \tilde q_{\overline A}(n).
\]
\end{proposition}
\begin{proof}
Same proof as Theorem~\ref{thm:Hilbert-poly}.
%Usual proof (by induction on $k$).
\end{proof}
We call $\tilde q_{\overline A}$ the intrinsic Hilbert polynomial of $\overline A$,
and denote by $\tildemu(\overline A)$ its leading term.

\begin{remark}
Assume that $\overline A$ is a $\lambda$-inert filtering on $A$.
Assume moreover that $\lambda(A_{0}) < \infty$.
In this situation, we  have defined the Hilbert polynomial $q_{\overline A}$.
We have, for every $n \in \N$
\[
\lambda(\tilde A_{n}) = \lambda(A_{n+1}) - \lambda(A_{n}),
\]
and therefore
\[
\tilde q_{\overline A} = \Delta q_{\overline A},
\]
where $\Delta p$ is the difference  of $p$: the polynomial defined by
$\Delta p(t) = p(t+1) - p (t)$.
\end{remark}

The intrinsic Hilbert polynomial becomes interesting when $\lambda(A_{0})$ is
infinite (and therefore we cannot compute the usual Hilbert polynomial).

\begin{definition}
Let $A_{0} \submod A$ be an $R$-submodule.  Denote
$\tilde A_{0} \coloneqq (S_{1} A_{0})/A_{0}$.  We say that $A_{0}$ is
\intro{$\lambda$-inert} if:
\begin{enumerate}
\item $\tilde A_{0}$ is finitely generated (as $R$-module);
\item $\lambda(\tilde A_{0}) < \infty$.
\end{enumerate}
\end{definition}
Assume that $A_{0}$ is $\lambda$-inert.
%\begin{enumerate}
%\item $A_{0}$ is $\lambda$-inert;
%\item $S A_{0} = A$;
%\end{enumerate}
We can define the associated filtering  $\filt(A_{0};A)$ of~$A$.
We denote
\[
\btAz \coloneqq \blowtilde (\filt(A_{0};A)) = \bigoplus_{n} (S_{n+1}A_{0})/(S_{n} A_{0})\, t^{n}.
\]
%, where $A_{n} \coloneqq  S_{n} A_{0}$.
%Then, $\Rees(\overline A)$ is Noetherian (since it is generated, as $S[y]$-module,
%by $A_{0} y^{0}$), and therefore $\overline A$ is a $\lambda$-inert filtration.
Then, $\btAz$ is a Noetherian $S$-module (since it is generated
by $\tilde A_{0}$), and therefore $\filt(A_{0},A)$ is a $\lambda$-inert filtering.

%satisfies the assumptions
%of Proposition~\ref{prop:intrinsic-Hilbert}.
Therefore, $\tilde q_{\btAz}$ exists; we define
\[
\tilde q_{A_{0}} \coloneqq q_{\btAz}, \qquad
 \tildemu[A_{0}] \coloneqq \mu(\btAz).
\]
%and $\tildemu[V_{0}]$ to be the leading term of $\tilde q_{V_{0}}$.

Thus, by definition, for $n$ large enough, 
\[
\tilde q_{A_{0}}(n) = \lambda(S_{n+1} A_{0}/S_{n} A_{0}).
\]

Unlike $\mu$, it can happen that  $\tildemu[A_{0}]$  depends on the
choice of $A_{0}$ (even when $S A_{0} = A$), as the following example shows (suggested by S. Virili):
\begin{example}\label{ex:Virili}
Let $R \coloneqq \Z$, $A \coloneqq S$, $\lambda$ be the standard length. 
%in Example~\ref{ex:length}(\ref{ex:length-Z}).
Fix $n \in \N$ and let 
\[
V \coloneqq \Z  + n S = \Z \oplus n\Z \x \oplus n\Z \x^{2} \oplus \cdots \submod A.
\]
Then, 
\[
\widetilde{V} = (V + \x V)/V \simeq \x \Z/n.
\]
Thus, %\lambda(\widetilde{V(n)}) 
$\lambda(\tilde V) = k\log n < \infty$; moreover, $V$ is finitely
generated and $S V = A$.
Moreover,
\[
V_{i} \coloneqq S_{i} V = S_{i} \Z + n S_{i} S = S_{i} + n S,
\]
and therefore
\[
\tilde V_{i} = (S_{i+1} + nS)/(S_{i} + nS) \simeq S^{(i + 1)}/(n).
\]
Thus,
\[
\lambda(\tilde V_{i}) = \log n \cdot \card{S^{(i+1)}} = \tbinom{i+k}{k-1} \log n 
\]
and therefore
\[
\tildemu[V] = \frac{\log n}{(k-1)!} \cdot t^{k-1}.
\]
Therefore,
\[
\tildemu(A) = \infty t^{k-1},
\]
and $\tildemu[V]$ depends in this case on the choice of~$V$.
\end{example}

\begin{definition}
%Assume that $A$ is a finitely generated $S$-module.
Define
\[
\tildemu(A) \coloneqq \sup \set{ \tildemu[ V_{0}]:\ 
V_{0} \submod A \text{ is $\lambda$-inert} 
%\ \ %\tilde V_{0} \text{ is finitely generated,  and}
%\ \ S V_{0} = A
} \in \vmc.
\]
\end{definition}
% (We will see later that we can drop the  assumption that
% $S V_{0} = A$ is the definition of $\tildemu(A)$).

\begin{lemma}\label{lem:inert-tilde}
Let $A$ be an $S$-module.
%Assume that $A$ is a finitely generated $S$-module.
Let $\overline A \coloneqq (A_{i})_{i \in \N}$ be a $\lambda$-inert filtering on~$A$.
Then, $\btA$ is an acceptable graded $S$-module
(see Definition~\ref{def:acceptable-grading})
%$\overline A$ is an acceptable filtration (see Definition~\ref{def:good-filter}), 
and therefore, 
%by Proposition~\ref{prop:AR-up}, 
by Proposition~\ref{prop:AR-graded-up},
there exists $d_{0} \in \N$ such that, for every
$n \in \N$, 
\[
\tilde A_{n+d_{0}} = S^{(n)} \tilde A_{d_{0}} = S_{n} A_{d_{0}}.
\]
Let $d_{1} \geq d_{0}$ such that $\lambda(\tilde A_{d_{1}}) < \infty$.
Then, $A_{d_{1}}$ is $\lambda$-inert %, $\tilde A_{d_{1}}$ is finitely generated,
and $S A_{d_{1}} = S A_{d_{0}}$.
Moreover,
\begin{equation}\label{eq:tilde}
\tildemu (\overline A) = \tildemu[A_{d_{1}}].
\end{equation}
Therefore, 
\[
\tildemu(A) = \sup \set{\tildemu[\overline A]: \ \overline A \text{ $\lambda$-inert
    filtering on } A},
\]
and 
if $A$ is finitely generated, then
\[
\tildemu(A) = \sup \set{\tildemu[\overline A]: \ \overline A \text{ $\lambda$-inert exhaustive
    filtering on } A}.
\]
\end{lemma}
\begin{proof}
It suffices to prove \eqref{eq:tilde}: the rest of the lemma is clear.
However, we have, for every $n \in \N$
\[
\tilde A_{n+d_{1}} = S_{n} \cdot  \tilde A_{d_{1}};
\]
 therefore,
\[
 \lambda(\tilde A_{n+d_{1}}) = \lambda(S_{n} \cdot  \tilde A_{d_{1}})
\]
proving that
\[
\tilde q_{\overline A}(n + d_{1}) = \tilde q_{A_{d_{1}}}(n)
\]
and thus \eqref{eq:tilde} follows.
\end{proof}

% \begin{definition}
% Let $\overline A$ be a filtering on $A$.
% Given $n \in \N$, the translate of $\overline A$ by $n$ is
% \[
% \overline A[n] \coloneqq \Pa{A_{n+i}}_{i \in \N}.
% \]
% \end{definition}

\begin{conjecture}
Let
$0 \longto  A \longto B \overset{\pi}\longto C \longto 0$ be an exact
sequence of %finitely generated 
$S$-modules.
Then,
\[
\tildemu(B) = \tildemu(A) \oplus \tildemu(C).
\]
\end{conjecture}
A particular case of the above conjecture (when $k = 1$ and $R = \Z$) is known
from \cites{DGSV:15, SV:18}.

We can prove sub-additivity quite easily:

\begin{thm}
Let
$0 \longto  A \longto B \overset{\pi}\longto C \longto 0$ be an exact
sequence of %finitely generated 
$S$-modules.
Then,
\[
\tildemu(B) \leq \tildemu(A) \oplus \tildemu(C).
\]
Moreover, $\tildemu(A) \leq \tildemu(B)$ and $\tildemu(C) \leq \tildemu(B)$.
\end{thm}
\begin{proof}
Let $\overline B = (B_{n})_{n \in \N}$ be a $\lambda$-inert filtering on $B$.
For every $n \in \N$, define 
\[\begin{aligned}
A_{n} &\coloneqq A \cap B_{n},  &
C_{n} &\coloneqq \pi(B_{n});\\
\overline A& \coloneqq (A_{n})_{n \in \N}, &
\overline C& \coloneqq (C_{n})_{n \in \N} = \filt(C_{0};C).
\end{aligned}\]

%Let $B_{0} \submod B$ be a $\lambda$-inert %finitely generated 
%$R$-submodule of $B$. 
% For every $n \in \N$, define 
% \[\begin{aligned}
% B_{n} &\coloneqq S_{n} B_{0}, &
% A_{n} &\coloneqq A \cap B_{n},  &
% C_{n} &\coloneqq \pi(B_{n});\\
% \overline B& \coloneqq (B_{n})_{n \in \N} = \filt(B_{0}; B), &
% \overline A& \coloneqq (A_{n})_{n \in \N}, &
% \overline C& \coloneqq (C_{n})_{n \in \N} = \filt(C_{0};C).
% \end{aligned}
% \]
It is clear that $\overline A$, $\overline B$, and $\overline C$ are  filterings on  
$A$, $B$, $C$, respectively.
\begin{claim}
$\overline C$ and $\overline A$ are $\lambda$-inert.
\end{claim}
In fact, for every $n \in \N$, we have an exact sequence of $R$-modules
\[
0 \to \tilde A_{n} \to \tilde B_{n} \to \tilde C_{n} \to 0.
\] 
Since $\tilde B_{0}$ is finitely generated and $R$ is Noetherian, both
$\tilde A_{0}$ and $\tilde C_{0}$ are finitely generated.
Since moreover (for $n$ large enough) $\lambda(\tilde B_{n})$ is finite, both 
$\lambda(\tilde A_{n})$ and $\lambda(\tilde C_{n})$ are finite.

Finally, we have the following  exact sequence of $S[y]$-modules:
\[
0 \to \blowtilde(\overline A) \to  \blowtilde(\overline B) \to \blowtilde (\overline C) \to 0;
\]
since $\blowtilde(\overline B)$ is Noetherian, also
$\blowtilde(\overline A)$ and $\blowtilde(\overline C)$ are Noetherian.

Notice moreover that, for every $n \in \N$ large enough,
\[
\tilde q_{\overline A }(n) + \tilde q_{\overline C }(n) =
\lambda(\tilde A_{n}) +  \lambda(\tilde C_{n}) 
= \lambda(\tilde B_{n}) = \tilde q_{\overline B}(n).
\]
Thus, 
\[
\tildemu[\overline B] = \tildemu(\overline A) \oplus \tildemu(\overline C).
\]
% \[
% \tildemu[B_{0}] = \tildemu(\overline A) \oplus \tildemu(\overline C) \le
% \tildemu(A) + \tildemu(C).
% \]
Therefore, by Lemma~\ref{lem:inert-tilde},
\[
\tildemu(B) \le
\tildemu(A) \oplus \tildemu(C). 
\]

\smallskip

Similarly, if $\overline A = (A_{n})_{n \in \N}$ is a $\lambda$-inert filtering on~$A$, 
then it is also a $\lambda$-inert filtering on~$B$, and therefore
$\tildemu[\overline A] \leq \tildemu(B)$; thus; $\tildemu(A) \leq \tildemu(B)$.

Finally, if $\overline C = (C_{n})_{n \in \N}$ is a $\lambda$-inert filtering on $C$, then
$\Pa{\pi^{-1}(C_{n})}_{n \in \N}$ is a
$\lambda$-inert filtering on~$B$, and therefore
$\tildemu[\overline C] \leq \tildemu(B)$; thus; $\tildemu(C) \leq \tildemu(B)$.
\end{proof}

Let $A$ be an $S$-module.
Let $d \in \N$ be the degree of $\tildemu(A)$ and $s \in \R$ be its coefficient.
We define the \intro{intrinsic \ldim} of $A$ to be
$d+1$ if $\tildemu(A) \ne 0$,
$0$ if $\lambda(A)  > 0$ and $\tildemu(A) = 0$, and   $-\infty$ if $\lambda(A) =0$.
For each $i \le k$, the \intro{intrinsic $i$-dimensional $\lambda$-entropy} of $A$ is
\[
\tilde h^{(i)}_{\lambda}(A) \coloneqq
\begin{cases}
0 & \text{if } i > d +1;\\
0 &\text{if } i \geq 1 \text{ and } \tildemu(A) = 0;\\
\infty & \text{if } i \leq d;\\
\frac{s}{(d+1)!} & \text{if } i = d+1 \text{ and } \tildemu(A) \ne 0;\\
\lambda(A) & \text{if } i = 0.
\end{cases}
\]
The intrinsic  $\lambda$-entropy is
\[
\tilde h_{\lambda}(A) \coloneqq \tilde h^{(k)}_{\lambda}(A)
\]
and has been studied already (at least, in the case $R = \Z$ and $k=1$)
in \cites{DGSV:15, GS:15, SV:18}, and for some non-Noetherian rings in \cite{SV:19}):
our definition of $\tildemu$ is clearly inspired by the intrinsic
algebraic entropy (the notation $\widetilde{\ent}$ is used elsewhere, but we
prefer $\tilde h_{\lambda}$ for consistency).

% We have therefore 2 length functions on $S$-modules extending $h_{\lambda}$ on locally \lR-finite
% $S$-modules: 
In general, $\hat h_{\lambda}$ and $\tilde h_{\lambda}$ are different.
\begin{example}
Let $R = \Z$, $k = 1$, $\lambda$ be the standard length,
and $A = \Q$ seen as $S$-module via the action
$x q \coloneqq q/2$.
Then,
$\tilde h_{\lambda}(A) = \log 2$ (see \cite{DGSV:15}),
% : but the computation is easy,
% since any finitely generated subgroup of $\Q$ is either $0$ or isomorphic to $\Z$),
while $\hat h_{\lambda}(A) = 0$
(notice that $\muhat(A) = \infty \cdot t^{0}$).
\end{example}

\section{Fine grading}\label{sec:fine}
Up to now, we have only considered the case when the degrees are natural
numbers.
As in the classical case when $R$ is a field, one can consider gradings in any
commutative monoid (see e.g. \cite{MS:05}).

\subsection{Graded modules}
Let $\Gamma$ be a commutative monoid.

% \begin{definition}
% We say that $\Gamma$ is ``good'' if:
% \begin{enumerate}
% \item $\Gamma$ is finitely generated (as a monoid);
% \item $\le$ is a partial ordering (\ie, if $m + p = n$ and $n + q = m$, then
% $m = n$);
% \item $\le$ is well-founded (CHECK).
% \end{enumerate}
% \end{definition}

Remember that $\Gamma$ has a canonical quasi-ordering, given by
$m \le n$ if there exists $p \in \Gamma$ with $m + p = n$.
The neutral element $0$ is a minimum of $\tuple{\Gamma, \leq}$.

\begin{definition}
We say that $\bar \gamma \in \Gamma^{k}$ is \textbf{good} (inside $\Gamma$) if:\\
for every 
%$k \in \N$,  $\bar n \in \N^{k}$, and 
$\lambda \in \Gamma$
there exist at most finitely many $\bar n  \in \N^{k}$, such that
$\bar n \cdot \bar \gamma = \lambda$.

We say that $\bar \gamma \in \Gamma^{k}$ is \textbf{very good} if:\\
for every 
%$k \in \N$,  $\bar n \in \N^{k}$, and 
$\lambda \in \Gamma$
there exist at most finitely many $\bar n  \in \N^{k}$, such that
$\bar n \cdot \bar \gamma \le \lambda$.
\end{definition}

For example, $\bar \gamma \in \Z^{k}$ is good (in $\Z$) if $\gamma_{i} > 0$ for $i = 1, \dotsc, k$.
$\bar \gamma \in \N^{k}$ is very good (in $\N$) iff $\gamma_{i} \ne 0$ for $i = 1, \dotsc, k$.
Notice that, in general, if $\bar \gamma$ is good, then each $\gamma_{i}$ is non-zero
(and even non-torsion).
%The canonical example of a good commutative monoid is $\N^{k}$, for $k \in \N$.

\begin{remark}
Given  $\bar \gamma = \tuple{\gamma_{1}, \dotsc, \gamma_{k}}\in \Gamma^{k}$, 
%a good commutative monoid $\Gamma$, and $\gamma_{1}, \dotsc, \gamma_{k} \in \Gamma$,
if $\bar \gamma $ is good, then 
the following expression is well defined:
\[
\frac{1}{\prod_{i=1}^{k}(1 - t^{\gamma_{i}})} \in \N[[t^{\Gamma}]]
\]
\end{remark}

\begin{definition}
Fix $\bar \gamma = \tuple{\gamma_{1}, \dotsc, \gamma_{k}} \in \Gamma^{k}$.
A $\Gamma$-graded $S$-module of degree $\bar \gamma$ is given by an $S$-module $M$ and a
decomposition 
\[
M = \bigoplus_{n\in \Gamma} M_{n},
\]
where each $M_{n}$ is an $R$-module, and, for every $i \le k$ and $n \in \Gamma$,
\[
x_{i} M_{n} \submod M_{n+ \gamma_{i}}.
\]
We denote by
$\overline M$ the module $M$ with the given grading (including the tuple  $\bar \gamma :=
\tuple{\gamma_{1}, \dotsc, \gamma_{k}}$).
\end{definition}

\begin{thm}\label{thm:Hilbert-graded-monoid}
Let $\overline M$ be a $\Gamma$-graded $S$-module of degree
$\bar \gamma \in \Gamma^{k}$.
For every $n \in \N$, let $a_{n} := \lambda(M_{n})$.
Define
\[
F_{\overline M}(t) := \sum_{n \in \Gamma} a_{n} t^{n}.
\]

Assume that:  
\begin{enumerate}
\item $\bar \gamma$ is good;
\item $\lambda(M_{n}) < \infty$ for every $n \in \Gamma$;
\item  $M$ is a Noetherian $S$-module.
\end{enumerate}
%and it is \lfgS.
%If $M$ is $\lambda$-Noetherian, 
Then, there exists a polynomial
$p(t) \in \R[t^{\Gamma}] $ such that
\[
F_{\overline M}(t) = \frac{p(t)}{\prod_{i=1}^{k}(1 - t^{\gamma_{i}})}.
\]
\end{thm}
\begin{proof}
\emph{Mutatis mutandis}, same proof as  Thm.~\ref{thm:Hilbert-graded}.
\end{proof}

\subsection{Filtered modules}

We move now from graded modules to  filtered modules.
%We need some condition on the grading monoid~$\Gamma$.

\begin{definition}
Let $\bar \gamma := \tuple {\gamma_{1}, \dotsc, \gamma_{k}} \in \Gamma^{k}$ and $N$ be an $S$-module.

An (increasing) \intro{$\Gamma$-filtering} on~$N$ with degrees~$\bar \gamma$ is a sequence
of $R$-submodules of $N$
\[
(N_{i}: i \in \Gamma)
\]
%every $N_{i}$ is contained in $N_{i+1}$ and 
such that it is increasing (\ie, if $i \le j$, then $N_{i} \submod N_{j}$),
$\bigcup_{i \in \Gamma} N_{i} = N$,
and $x_{i} N_{j} \submod N_{j + \alpha_{i}}$ for every $j \in \Gamma$, $i \le k$.
We denote by $\overline N$ the $S$-module with the given tuple $\gamma$ and the
filtering
$(N_{i})_{i \in \Gamma}$.
\end{definition}

\begin{definition}
Let $\bar \delta = \tuple{\delta_{1}, \delta_{2}, \dots}$ be a tuple of generators of~$\Gamma$
(for simplicity, we assume $\Gamma$ countable: later we will be interested only in
the case when $\bar \delta$ is a finite tuple).
Let $\bar y$ be a tuple of variables indexed by~$\bar \delta$
(\ie, there is one variable $y_{j}$ for each chosen generator $\delta_{j}$).
The \intro{blow-up} module associated to $\overline N$ and the tuple $\bar \delta$ is
the following graded $S[\y]$-module.
As an $R$-module, 
\[
\RNd := \bigoplus_{n \in \Gamma} N_{n} t^{n}.
\]
The $\Gamma$-grading of $\RNd$ is given by the decomposition 
$\RNd = \bigoplus_{n \in \Gamma} N_{n} t^{n}$.

The multiplication by $x_{i}$ on $\RNd$ is defined as:
\[
x_{i}(v t^{n}) := (x_{i}v) t^{n + \gamma_{i}},
\]
for every $i \le k$, $n \in \Gamma$, $v \in N_{n}$, and then extended by $R$-linearity
on all $\RNd$: notice that the  $x_{i}$  has degree $\gamma_{i}$ on~$\RNd$.
For each $\delta_{j} \in \bar \delta$, the multiplication by $y_{j}$ on $\RNd$ is defined as: 
\[
y_{j} (v t^{n}) := v t^{n+\delta_{j}},
\]
for every $n \in \Gamma$, $v \in N_{n}$, and then extended by $R$-linearity
on all $\RNd$: notice that $y_j$ has degree~$\delta_{j}$.
\end{definition}

\begin{thm}\label{thm:Hilbert-filtered-monoid}
Let $\bar \N$
%$N = \bigcup N_{i}$ 
be a $\Gamma$-filtering on $N$ with degrees $\bar \gamma$.
%denoted by $\overline N$.
%Assume that $\Gamma$ is good, and let 
Let $\bar \delta = \tuple{\delta_{1}, \dotsc, \delta_{\ell}}$ be a
finite tuple of generators of $\Gamma$, with corresponding variables
$\bar y = \tuple{y_{1}, \dotsc, y_{\ell}}$.

For every $n \in \N$, let $a_{n} := \lambda(N_{n})$.
Define
\[
F_{\overline N}(t) := \sum_{n \in \Gamma} a_{n} t^{n} \in \R[[t^{\Gamma}]]
\]

Then,
\[
F_{\overline N}  = F_{\RNd}.
\]

Therefore, if we assume that:
\begin{enumerate}
\item $\bar \gamma \cup \bar \delta$ is good;
\item for every $n \in \N$, $\lambda(N_{n}) < \infty$;
\item $\ReesN$ is Noetherian as $S[\bar y]$-module.
\end{enumerate}
Then, there exists a polynomial $p(t) \in \R[t^{\Gamma}]$ such that
\[
F_{\overline N}(t) = \frac{p(t)}{ \prod_{j = 1}^{\ell}(1 - t^{\delta_{j}}) \prod_{i = 1}^{k} (1 - t^{\gamma_{i}})}
\]
(the $(1 - t^{\delta_{j}})$-factor in the denominator is due to the action of
$y_{j}$ on $\RNd$ of degree~$\delta_{k}$).
\end{thm}

\subsection{Growth function}
In this subsection we fix a monoid $\Gamma$ with a tuple of generators
$\bar \delta = \tuple{\delta_{1}, \dotsc, \delta_{\ell}}$.

We also fix a tuple $\bar \gamma = \tuple{\gamma_{1}, \dotsc, \gamma_{k}} \in \Gamma^{k}$.
Given a monomial  in $S = R[x_1, \dotsc, x_{k}]$ its $\bar \gamma$-degree 
$\deg_{\bar\gamma}$ is defined in the ``obvious'' way:
\[
\deg_{\bar \gamma}(r x_{1}^{n_{1}} \cdots x_{k}^{n_{k}}) :=
n_{1} \gamma_{1} + \dots + n_{k}\gamma_{k}.
\]
Given a polynomial $p(\x) \in S = R[x_1, \dotsc, x_{k}]$, we say that 
its $\bar \gamma$-degree is less or equal to $n \in \Gamma$, and write $\deg_{\bar \gamma}(p) \le n$, if
each monomial in $p$ has $\bar \gamma$-degree less or equal to~$n$
(since $\Gamma$ is not linearly ordered in general, it's not clear how to define
the $\bar \gamma$-degree of a polynomial).
For every $n \in \Gamma$, we denote
\[
S_{n} := \tuple{p \in S: \deg_{\gamma}(p) \le n}.
\]

\begin{thm}\label{thm:Hilbert-orbit-monoid}
Let $N$ be an $S$-module.
Let $V_{0} \submod N$ be an $R$-submodule.
For every $n \in \Gamma$, let $V_{n} := S_{n} V_{0}$ (notice that $S_{0} =R$, that
$V_{0} = S_{0} V_{0}$, and that $S_{n}$ and $V_n$ are $R$-modules), and $a_{n} := \lambda(V_{n})$.
Define
\[
G_{V_{0}}(t) := \sum_{n \in \Gamma} a_{n} t^{n} \in \R[[t^{\Gamma}]] 
\]
Assume that:
\begin{enumerate}
\item $\bar \gamma \cup \bar \delta$ is very good (inside $\Gamma$);
\item $\lambda(V_{0}) < \infty$;
\item $V_{0}$ is finitely generated as $R$-module.
%\item $R$ is a Noetherian ring.
\end{enumerate}
Then, each $a_{n}$ is finite, and there exists a polynomial
$p(t) \in \R[t^{\Gamma}]$ such that
\[
G_{V_{0}}(t) = \frac{p(t)}{\prod_{j = 1}^{\ell} (1 - t^{\delta_{j}}) \prod_{i = 1}^{k} (1 - t^{\gamma_{i}})}
\]
\end{thm}
\begin{proof}
The fact that $\gamma$ is very good is equivalent to the fact that $S_{n}$ is
finitely generated (as $R$-module) for every $n \in \Gamma$.
The above plus the fact that $\lambda(V_{0}) < \infty$ easily implies that 
$\lambda(V_{n}) < \infty$ for every $n \in \Gamma$.
Let $\overline V := S V_{0}$ as filtered $S[\y]$-module.
Then, $V_{0} t^{0}$ generates $\Rees(\overline V)$ as $S[\y]$-module, and
therefore
$\Rees(\overline V)$ is a Noetherian $S[\y]$-module.

We can conclude as in the proof of Theorem~\ref {thm:Hilbert-orbit}, using 
Theorem~\ref{thm:Hilbert-filtered-monoid}.
\end{proof}

\begin{example}
Let $\x := \tuple{x_{1}, \dotsc, x_{k}}$, $\y := \tuple{y_{1}, \dotsc, y_{\ell}}$,
$S := R[\x,\y]$, $\Gamma := \N^{2}$,
$\gamma_{i} := \delta_{1} := \tuple{1,0}$ for $i = 1, \dotsc, k$, and $\gamma_{i} := \delta_{2} :=
\tuple{0,1}$ for
$i := k+1, \dotsc, k + \ell$.
Thus, each $x_{i}$ has degree $\delta_{1}$ and each $y_{j}$ has degree
$\delta_{2}$.
A monomial in $\x\y$ has therefore a ``double degree'' $\tuple{m,n} \in \Gamma$, where
$m$ is its total degree in $\x$ and $n$ is its total degree in $\y$.
A polynomial in $\R[t^{\Gamma}]$ is the same object as a polynomial in the two
variables
$t_{1}, t_{2}$.
Let $N$ be an $S$-module and $V_{0}\submod N$ be an $R$-submodule which satisfies
(2) and (3).
Then, if $R$ is Noetherian, we have a corresponding function
\[
G_{V_{0}}(t_{1}, t_{2}) = \frac{p(t_{1}, t_{2})}
{(1-t_{1})^{k+1} (1 - t_{2})^{\ell+1}}
\]
where $p \in \R[t_{1}, t_{2}]$.
\end{example}

\subsection{Multi-variate Hilbert polynomial}
Let $P \coloneqq \tuple{P_{1}, \dotsc, P_{\ell}}$ be a partition of 
$\set{1, \dotsc, k}$ into $\ell$ nonempty subsets;
for every $j \le \ell$, let $p_{j}$ be the cardinality of $P_{j}$.
In the following, we will assume that
$P_{1} = \set{1, 2, \dotsc, p_{1}}$,
$P_{2} = \set{p_{1}+1, p_{1}+2,\dotsc, p_{1} + p_{2}}$,
\dots,
$P_{\ell} = \set{p_{1}+ \dots + p_{\ell-1} + 1, \dotsc, k}$.

Let $\Gamma := \N^{\ell}$;
for every $j \leq \ell$, let $\hat e_{j} \in \Gamma$ be the vector with all
coordinates $0$ except the $j$-th which is $1$.
Let $\x := \tuple{x_{1}, \dotsc, x_{k}}$ be a $k$-tuple of variables;
to each variable $x_{i} \in P_{j}$ assign the weight $\hat e_{j}$,
and define
\[
\bar e \coloneqq \tuple{\hat e_{1}, \hat e_{1}, \dotsc, \hat e_{\ell}} \in \N^{k},
\]
where each weight $\hat e_{j}$ is repeated $p_{j}$ times.

As usual,  $S \coloneqq  R[x_{1}, \dotsc, x_{k}]$;
for every $\bar m \in \N^{\ell}$, let 
\[
S^{(\bar e)}_{\bar m} := \set{p \in S: \deg_{\bar e}(p) \le \bar m}.
\]
An equivalent way of describing $S_{\bar m}^{(\bar e)}$ is the following.
Let $\bar t := \tuple{t_{1}, \dotsc, t_{\ell}}$.
Let $\phi: S \to R[\bar t]$ be the homomorphism of $R$-algebrae mapping
$x_{i}$ to $t_{j}$ when $i \in P_{j}$.
Then, $q \in S_{\bar m}^{(\bar e)}$ iff, for every $j \le \ell$,
$\deg_{t_{j}}(\phi(q)) \le m_{j}$.

Let $M$ be a module over $S$ and
 $V_{0} \submod M$ be an $R$-submodule.
For every $\bar m \in \N^{\ell}$, define $V_{\bar m} := S_{\bar m}^{\bar e}V_{0}\submod  M$ 
and
$a_{\bar m} := \lambda(V_{\bar m})$.
\begin{thm}
In the above setting, assume that $V_{0} \submod M$ witnesses that
$M$ is $\lambda_{S}$-finite.
% \begin{enumerate}
% \item $\lambda(V_{0}) < \infty$;
% \item $V_{0}$ is finitely generated as $R$-module;
% \item $R$ is a Noetherian ring;
% \item $S V_{0} = M$
% \end{enumerate}
Then, each $a_{\bar m}$ is finite.
Moreover, there exists a polynomial
$q(\bar t) \in \R[\bar t]$ such that: 
\begin{enumerate}[(i)]
\item for every $\bar m \in \N^{\ell}$ large enough,
\[
a_{\bar m} = q(\bar m);
\]
\item for every $j \le \ell$,
\[
\deg_{t_{j}}(q) \le p_{j}.
\]
\end{enumerate}
Moreover, the leading homogeneous component of $q$ (see Def.~\ref{def:leader})
is independent from~$V_{0}$.
\end{thm}
\begin{proof}
Choose $\bar e$ as tuple of generators of $\Gamma$ and apply
Theorem~\ref{thm:Hilbert-orbit-monoid}.
We obtain that there exists a polynomial $r(\bar t) \in \R[\bar t]$ such that
\[
\sum_{\bar m \in \N^{\ell}} a_{\bar m} \bar t^{\bar m} =
\frac{p(\bar t)}{\prod_{j = 1}^{\ell} (1 - t_{j})^{p_{j} + 1}} 
\]
(the exponents $p_{j} + 1$ in the denominator come from the combination
of the $p_{j}$ variables in $P_{j}$, each with degree $\hat e_{j}$,
plus the generator $\hat e_{j}$).
%Any expression as above satisfies the conclusion of the theorem:
By Proposition~\ref{prop:poly-growth},
there exists a polynomial $q$ satisfying (i) and (ii).
By Proposition~\ref{prop:leader}, the leading homogeneous component of $q$ is
independent from~$V_{0}$.
\end{proof}

%\todo{Give application with $\Gamma = \N^{k}$, and some references.}

\appendix
%\clearpage

\section{Main definitions}
For the reader convenience, we  collect here the various definitions used up 
to~\S\ref{sec:intrinsic}.
\begin{adef}{Ring}
$R$ is a commutative ring with $1$ (most of the time, Noetherian).\\
$S = R[x_{1}, \dotsc, x_{k}]$ for a fixed $1 \leq k \in \N$.\\
$S_{n}$ is the set of polynomials of degree at most $n$, and
$S^{(n)}$ is the set of homogeneous polynomials of degree exactly $n$ (plus $0$).
\end{adef}

\subsection{Lenght}
\begin{adef}{Lenght}
A length  is a function
\[
\lambda: \Rmod \to \Rgeqzinf
\]
satisfying the following conditions:
\begin{enumerate}
\item $\lambda(0) = 0$;
\item $\lambda(M) = \lambda(M')$ when $M$ and $M'$ are isomorphic;
\item for every exact sequence
$0 \to A \to B \to C \to 0$,
\[
\lambda(B) = \lambda(A) + \lambda(C);
\]
\item for every $M \in \Rmod$,
\[
\lambda(M) = \sup \set{\lambda(M'): M' \submod M\text{ finitely generated $R$-submodule}}.
\]
\end{enumerate}
\end{adef}
Fix a length function $\lambda$ and an $R$-module $N$.
\begin{adef}{Locally $\lR$-finite}
$N$ is locally $\lR$-finite if, for every $N' \submod N$ finitely
generated, $\lambda(N')$ is finite.
\end{adef}

Let $T$ be an $R$-algebra and $A$ be a $T$-module.
\begin{adef}{$\lT$-good}
$A$ is \lT-small if it is locally $\lR$-finite and it is finitely generated as
$T$-module; a witness for it is a finitely generated $R$-submodule $V$ such
that $\lambda(V) < \infty$ and $TV = A$.
\end{adef}

\subsection{Graded  modules}
Fix $\bar \gamma = \tuple{\gamma_{1}, \dotsc, \gamma_{k}} \in \N^{k}$ and an $S$-module $M$.
\begin{adef}{Grading}
An $\N$-grading of $M$  of degree $\bar \gamma$ is given by a
decomposition 
\[
M = \bigoplus_{n \in \N} M_{n},
\]
where each $M_{n}$ is an $R$-module, and, for every $i \le k$ and $n \in \N$,
\[
x_{i} M_{n} \submod M_{n+ \gamma_{i}}.
\]
\end{adef}
\begin{adef}{Acceptable}
A graded $S$-module $\overline M$
is acceptable if $\gamma_{i} = 1$ for $i = 1, \dotsc, k$ and $M$ is finitely
generated as $S$-module.
\end{adef}

\subsection{Filtered module}
\begin{adef}{Filtering}
An (increasing) filtering $\overline M$ on $M$ with degrees~$\bar \gamma$ is an increasing sequence
of $R$\hyph{}submodules of $M$
\[
M_{0} \submod M_{1} \submod M_{2} \submod \dotsb \submod M
\]
%every $N_{i}$ is contained in $N_{i+1}$ and 
such that 
%$\bigcup_{i=0}^\infty M_{i} = M$, and 
$x_{i} M_{j} \submod M_{j + \gamma_{i}}$ for every $j \in \N$, $i \le k$, and it is
exhaustive if $\bigcup_{i=0}^\infty M_{i} = M$.
(Remember that all filterings we consider are increasing).
\end{adef}

Fix a filtering $\overline M$ on $M$.
% \begin{adef}{Translate}
% The translate of $\overline M$ by $n$ is the filtering
% \[
% \overline M[n] \coloneqq \Pa{M_{n+i}}_{i \in \N}.
% \]
% \end{adef}
\begin{adef}{Blow up}
The blow up of $\overline M$ is the graded $S[y]$-module of degree
$(\gamma_{1}, \dotsc, \gamma_{k}, 1)$
given by:
\[\Rees(M) = \bigoplus_{n \in \N} M_{n} y^{n},\] with scalar multiplication defined by
\[
x_{i}(v y^{j}) := (x_{i}v) y^{j + \gamma_{i}}, \qquad
y (v y^{j}) := v y^{j+1}.
\]
\end{adef}
\begin{adef}{Tightly generates}
Given $m \in \N$, define
\[
M^{m} := \bigoplus_{n \leq m} M_{n} y^{n} \submod \Rees(M).
\]
We say that $M_{m}$ tightly generates $\overline M$ if $M^{m}$ generates
$\Rees(\overline M)$.%
\footnote{By Lemma~\ref{lem:tight}, this is equivalent to the original definition}
\end{adef}

\begin{adef}{Induced filtration}
Given and $R$-submodule $V_{0} \submod M$, we have the filtering and the graded module
\begin{align*}
\filt(V_{0}; M) &\coloneqq (S_{n }V_{0})_{n \in \N}\\
\gr(V_{0}; M) &\coloneqq \Rees(\filt(V_{0}; M)) =
\bigoplus_{n \in \N} S_{n} V_{0}\, y^{n}.
\end{align*}
\end{adef}

\begin{adef}{Acceptable}
$\overline M$ is acceptable if  $\Rees(\overline M)$ is an acceptable graded
$S[y]$-module: that is, all variables $x_{i}$ and $y$ have degree $1$, 
and $\Rees(\overline M)$
is a finitely generated $S[y]$-module.

$\overline M$ is good if moreover $\forall n \in \N$ $\lambda(M_{n}) < \infty$.
% $\Rees(\overline M)$ is finitely generated (as an $S[y]$-module).
\end{adef}

\subsection{Value monoid}
\begin{adef}{Value monoid}
$\vmc$ is the set of monomials of the form $r t^{d}$, with
$r \in R_{>0} \cup \set \infty$ and $d \in \set{1, \dotsc, k}$ (plus $0$).
\end{adef}

\subsection{Hilbert  polynomial}
Let $V \submod M$ be an $R$-submodule of $M$.

\begin{adef}{Hilbert polynomial}
The Hilbert polynomial of $V$ is the polynomial $q_{V}(t) \in \R[t]$ (which exists
under suitable assumptions) such that, for every $n$ large enough,
\[
\lambda(S_{n} V) = q_{V}(n).
\]
\end{adef}

\begin{adef}{$\mu$}
If $V$ witnesses that $M$ is \lS-small, the leading term $\mu(V) = \mu_{\lambda}(V) \in \vm$ of
$q_{V}$ does not depend on the choice of $V$.
If $\mu(V) = m t^{d}$,
the $\lambda$-dimension of $M$ is $\lambdadim(M) \coloneqq d$, 
the $\lambda$-degree of $M$ is $d! m$.

If $M$ is not $\lambda_{S}$-small, 
\[
\mu(M) := \sup\set{\mu(M'): M' \text{ \lS-small submodule of } M} \in \vmc.
\]
\end{adef}
\begin{adef}{$d$-dimensional entropy}
For every $i \leq k$, the $i$-dimensional entropy of $M$ is
\[
h^{(i)}_{\lambda}(M) \coloneqq
\begin{cases}
\infty &\text{if } d > i\\
0 &\text{if } d < i\\
i! \, m &\text{if } d = i.
\end{cases}
\]
$h^{(1)}$ is the receptive entropy, $h^{(d)}$ is the algebraic entropy.
\end{adef}

\subsection{Hilbert-Samuel polynomial}
Let $I := (x_{1}, \dotsc, x_{k}) \ideal R$.
\begin{adef}{Hilbert_Samuel polynomial}
The Hilbert-Samuel polynomial of $M$ is the polynomial $\bar q(t) \in \R[t]$ (which exists
under suitable assumptions) such that, for every $n$ large enough,
\[
\lambda(M / I^{n+1} M) = \bar q(n).
\]
The leading term of $\bar q$ is $\bar \mu(M)$.
\end{adef}

\subsection{Modules over $R$-algebrae}
$T$ is a commutative $R$-algebra generated by $\bar\gamma \coloneqq \tuple{\gamma_{1},
  \dotsc, \gamma_{k}}$. 
Let $\phi: S \to T$ be the homomorphism of $R$-algebrae mapping $x_{i}$ to $\gamma_{i}$.
$A$ is a $T$-module.
\begin{adef}{$\mu(M;\phi)$}
We denote $A$ as an $S$-module by $\tuple{A;\bar \gamma} = \tuple{A; \bar \phi}$, and
correspondingly $\mu(A; \bar\gamma) \coloneqq \mu(\tuple{A;\bar \gamma})$ and
$\mu(A; \phi) \coloneqq \mu(\tuple{A;\phi})$.  
If $A$ is $\lT$-small
the dimension of $A$ (as $T$-module) is the degree of the monomial $\mu(A;
\bar\gamma)$, which does not depend on the choice of $\bar \gamma$.

For $i \leq k$, the $i$-dimensional entropy of $A$ \wrt $\bar\gamma$ is
\[
h^{(i,\bar\gamma)} \coloneqq \sup\set{h^{(d)}(A, \bar \gamma):
A' \submod A \text{ \lT-small $T$-submodule}
}.
\]
\end{adef}

\subsection{V\'amos construction}%{$\hat \mu$}
\begin{adef}{$\lS$-small chains}
An $\lS$-small chain in $M$ of size $n$ is a sequence of $S$-submodules
\[
\Mseq = \Pa{M_{1} \submod M_{2} \submod \dots \submod M_{2n-1} \submod M_{2n}
  \submod M},
\]
such that, for every $i =1, \dotsc, n$, $\lambda(\hat M_{i}) < \infty$, where
$\hat M_{i} \coloneqq M_{2 i}/M_{2 i -1}$.
\end{adef}
\begin{adef}{$hat \theta$}
Let $\theta$ be a suitable function.
\begin{align*}
\theta(\Mseq) &\coloneqq \sum_{i=1}^{n} \theta(\hat M_{i}),\\
\thetahat(A) &\coloneqq \sup\set{\theta(\Aseq): \text{ $\Mseq$ $\lS$-small chain in
    $M$}}.
\end{align*}
\end{adef}
\begin{adef}{$\hat \lambda$}
The corresponding quantities for $\mu$, $h_{\lambda}^{(i)}$ and
$h_{\lambda}^{(i,\bar \gamma)}$ are denoted by $\hat \mu_{\lambda}$,
$\hat h_{\lambda}^{(i)}$ and $\hat h_{\lambda}^{(i,\bar \gamma)}$, respectively.
\end{adef}

\subsection{Intrinsic entropy}
Let $\overline M$ be a filtering on $M$ and $V_{0} \submod M$.
\begin{adef}{Inert filtration}
\[\begin{aligned}
\tilde M_{i} &\coloneqq M_{i+1}/M_{i}\\
\blowtilde (\overline M) &\coloneqq \bigoplus_{i \in \N} {\tilde M_{i}}\\
\tilde{\gr}(V_{0}; M) &\coloneqq \blowtilde (\filt(V_{0};M)) = \bigoplus_{n} (S_{n+1}V_{0})/(S_{n} V_{0})\, y^{n}.
\end{aligned}\]
$\overline M$ is $\lambda$-inert if $\blowtilde (\overline M)$ is Noetherian and, for every $n$
large enough, $\lambda(\tilde M_{i}) < \infty$.\\
$V_0$ is $\lambda$-inert if $\lambda(\tilde V_{0}) < \infty$ and $\tilde V_{0}$ is finitely generated,
where $\tilde V_{0} \coloneqq (S_{1} V_{0}) / V_{0}$.
\end{adef}
\begin{adef}{Intrinsic Hilbert polynomial} 
The intrinsic Hilbert polynomial of the $\lambda$-inert filtering $\overline M$ is the
polynomial $q_{\overline M}(t) \in \R[t]$ (which exists
under suitable assumptions) such that, for every $n$ large enough,
\[
\lambda(\tilde M_{n}) = \tilde q_{\overline M}(n);
\]
$\tildemu(\overline M)$ is the leading term of $\tilde q_{\overline M}$.
\end{adef}
\begin{adef}{Intrinsic Hilbert polynomial} 
$\tilde q_{V_{0}} \coloneqq \tilde q_{\tilde{\gr}(M; V_{0}) }$ and 
$\tildemu[V_{0}]$ is the leading term of $\tilde q_{V_{0}}$.
\[
\tildemu(M) \coloneqq \sup \set{ \tildemu[ V_{0}]:\ 
V_{0} \submod M \text{ is $\lambda$-inert}}
\]
Let $d$ be the degree of $\tildemu(M)$ and $s$ be its coefficient.\\
The intrinsic $\lambda$-dimension of $M$ is 
\[
\begin{cases}
d+1 &\text{ if } \tildemu(M) \neq 0,\\ 
0  &\text{ if }  \lambda(M) >0  \text{ and } \tildemu(M) = 0\\
- \infty & \text{ if } \lambda(A) = 0.
\end{cases}
\]
The intrinsic $i$-dimensional $\lambda$-entropy of $M$ is
\[
\tilde h^{(i)}_{\lambda}(M) \coloneqq
\begin{cases}
0 & \text{if } i > d +1 \text{ or } (i = d + 1 \text{ and } \tildemu(M) = 0);\\
\infty & \text{if } i \leq d;\\
\frac{s}{(d+1)!} & \text{if } i = d+1 \text{ and } \tildemu(M) \ne 0;\\
\lambda(M) & \text{if } i = 0.
\end{cases}
\]
The intrinsic  $\lambda$-entropy is
$\tilde h_{\lambda} \coloneqq \tilde h^{(k)}_{\lambda}$.
\end{adef}
\subsection{Polynomials}
\begin{adef}{Homogeneous component}
Given $p(\bar t) \in \R[\bar t]$, the homogeneous component of $p$ of degree $i$
is denoted by $p_{i}$.  The leading homogeneous component of $p$ is $p_{d}$,
where $d$ is the degree of $p$ (or $0$ if $p = 0$).
\end{adef}

\section{Non-commuative rings}\label{sec:noncommutative-hat}
In this appendix, $R$ is no longer assumed to be commuative.
Let $G$ be some (associative) monoid, and
$T \coloneqq R[G]$. 
%is some $R$-algebra (not necessarily generated by central elements).
$\Tmod$ is the category of left $T$-modules (by ``modules'' we will mean left modules).
$\theta$ is some function from $\Tmod$ to %$\vmc$.
the family of monomials of the form
$r t^{n}$, with $n \in \N$ and $r \in \R_{> 0} \cup \set \infty$ (plus the monomial~$0$).
We assume:
\begin{enumerate}
\item $\theta$ is additive on the category of
\lT-small $T$-modules:
$\theta(0) = 0$ and, for every exact sequences of \lT-small $T$-modules
$0 \to A \to B \to C \to 0$, $\theta(B) = \theta(A) \oplus \theta(C)$;
\item $\theta$ is invariant: if $A$ and $B$ are isomorphic
$T$-modules, then $\theta(A) = \theta(B)$;
\item 
$
\theta(A) = \sup \set{\theta(B): B \submod A \text{ \lT-small $T$-submodule}}.
$
\end{enumerate}
We define
\[
\thetahat(A) \coloneqq \sup \set{\theta(\Aseq): \Aseq  \text{ \lT-small chain in $A$}}.
\]

\begin{thm}\label{thm:noncommutative-hat}
Assume that $T$ is (left) Noetherian.
Then,
\begin{enumerate}[(a)]
\item $\thetahat$ is invariant and additive on all $\Tmod$;
\item if $A$ is locally \lR-finite, then $\thetahat(A) = \theta(A)$;
\item $\thetahat(A) = \sup\set{\thetahat(B): B \submod A \text{ finitely 
generated $T$-submodule}}$;
\item $\thetahat \geq \theta$.
\end{enumerate}
\end{thm}
\begin{proof}
Same as Propositions~\ref{prop:theta} and~\ref{prop:cosmall}.
The assumption that $T = S[G]$ is used in the following way:
\begin{claim}
Let $A$ be a $T$-module, generated (as a $T$-module) by $a_{1}, \dotsc,
a_{\ell}$.
Assume that $\lambda(Ra_{i})$ is finite, for $i = 1, \dotsc, \ell$.
Then, $A$ is locally \lR-finite.
%$Ann_{R}(A) = Ann_{R}(a_{1}) \cap \dots \cap Ann_{R}(a_{\ell})$.
\end{claim}
\end{proof}
It is unclear if, under the same assumptions as in
Thm.~\ref{thm:noncommutative-hat}, we can conclude that, when $A$ is a
finitely generated $T$-module,
\[
\thetahat(A) = \sup \set{\thetahat(A/I): I \ideal R \text{ \lcosmall ideal}}.
\]

\begin{example}
Let $G$ be an amenable cancellative monoid and $T \coloneqq R[G]$.
Let $\lambda$ be a length on $\Rmod$.
Let $h_{\lambda}(A)$ be the algebraic entropy of the action of $G$ on $A$ for the
length  $\lambda$ (see \eg \cite{DFG}).
Assume that $T$ is Noetherian.
Then, $\theta \coloneqq h_{\lambda}$ satisfies the assumptions of this section.
Thus, the function $\hat h_{\lambda}$ satisfies the conclusion of
Theorem~\ref{thm:noncommutative-hat}, and in particular is a length function
on all $T$-mod.
\end{example}

%%%%%%%%%%%%%%%%%%%%%%%%%%%%%%%%%%%%%%%%%%%%%

%\bibliographystyle{alpha}	% uses file "filename.bst"
%\bibliography{hilbert}		% expects file "filename.bib"

\printbibliography

\end{document}